\documentclass[12pt]{article}
\usepackage{amssymb, amsmath, amsthm, amscd}
\usepackage[pdftex]{graphics}
\usepackage[utf8]{inputenc}
\usepackage[all,cmtip]{xy}
\usepackage{enumitem}

\usepackage[mathscr]{euscript}
\usepackage[colorlinks=true]{hyperref}
\hypersetup{colorlinks   = true}
\hypersetup{linkcolor=blue}
\usepackage{tikz}

\begin{document}

\newtheorem{The}{Theorem}[section]
\newtheorem{Lem}[The]{Lemma}
\newtheorem{Prop}[The]{Proposition}
\newtheorem{Cor}[The]{Corollary}
\newtheorem{Rem}[The]{Remark}
\newtheorem{Obs}[The]{Observation}
\newtheorem{SConj}[The]{Standard Conjecture}
\newtheorem{Titre}[The]{\!\!\!\! }
\newtheorem{Conj}[The]{Conjecture}
\newtheorem{Question}[The]{Question}
\newtheorem{Prob}[The]{Problem}
\newtheorem{Def}[The]{Definition}
\newtheorem{Not}[The]{Notation}
\newtheorem{Claim}[The]{Claim}
\newtheorem{Conc}[The]{Conclusion}
\newtheorem{Ex}[The]{Example}
\newtheorem{Fact}[The]{Fact}
\newtheorem{Formula}[The]{Formula}
\newtheorem{Formulae}[The]{Formulae}
\newtheorem{The-Def}[The]{Theorem and Definition}
\newtheorem{Prop-Def}[The]{Proposition and Definition}
\newtheorem{Lem-Def}[The]{Lemma and Definition}
\newtheorem{Cor-Def}[The]{Corollary and Definition}
\newtheorem{Conc-Def}[The]{Conclusion and Definition}
\newtheorem{Terminology}[The]{Note on terminology}
\newcommand{\C}{\mathbb{C}}
\newcommand{\R}{\mathbb{R}}
\newcommand{\N}{\mathbb{N}}
\newcommand{\Z}{\mathbb{Z}}
\newcommand{\Q}{\mathbb{Q}}
\newcommand{\Proj}{\mathbb{P}}
\newcommand{\Rc}{\mathcal{R}}
\newcommand{\Oc}{\mathcal{O}}
\newcommand{\Vc}{\mathcal{V}}
\newcommand{\Id}{\operatorname{Id}}
\newcommand{\pr}{\operatorname{pr}}
\newcommand{\rk}{\operatorname{rk}}
\newcommand{\del}{\partial}
\newcommand{\delbar}{\bar{\partial}}
\newcommand{\Cdot}{{\raisebox{-0.7ex}[0pt][0pt]{\scalebox{2.0}{$\cdot$}}}}
\newcommand\nilm{\Gamma\backslash G}
\newcommand\frg{{\mathfrak g}}
\newcommand{\fg}{\mathfrak g}
\newcommand{\Oh}{\mathcal{O}}
\newcommand{\Kur}{\operatorname{Kur}}
\newcommand\gc{\frg_\mathbb{C}}
\newcommand\slawek[1]{{\textcolor{red}{#1}}}
\newcommand\dan[1]{{\textcolor{blue}{#1}}}
\begin{center}
	
	{\Large\bf A Variational Approach to SKT and Balanced Metrics}
	
\end{center}

\begin{center}
	{\large S\l{}awomir Dinew and Dan Popovici}

\end{center}

\begin{center}

{\it To the memory of Jean-Pierre Demailly}

\end{center}

\vspace{1ex}

\noindent{\small{\bf Abstract.} We investigate compact complex manifolds endowed with SKT or balanced metrics. In each case, we define a new functional whose critical points are proved to be precisely the K\"ahler metrics, if any, on the manifold. As general manifolds of either type need not admit K\"ahler metrics, this provides an approach to new obstructions to K\"ahlerianity within these two families of metrics.}

\vspace{1ex}

\noindent{\small{\bf R\'esum\'e.} Nous \'etudions les vari\'et\'es complexes compactes lisses munies de m\'etriques hermitiennes sp\'eciales du type SKT ou \'equilibr\'e. Dans chacun de ces deux cas, nous d\'efinissons une fonctionnelle dont nous d\'emontrons que les points critiques, s'ils existent, sont exactement les m\'etriques k\"ahl\'eriennes, au cas o\`u elles existent, sur la vari\'et\'e. Cette approche constitue implicitement une recherche de nouvelles obstructions \`a l'existence de m\'etriques k\"ahl\'eriennes dans chacune de ces deux classes de  m\'etriques hermitiennes.}

\vspace{1ex}


\vspace{1ex}


\section{Introduction}

The study of complex manifolds admitting Hermitian metrics with special properties is currently one of the central topics in complex geometry. While the development of K\"ahler geometry is the premier example, the literature abounds in K\"ahler-like metric properties of complex manifolds for which a rich theory is being developed -- see, for example, [AB95], [FPS04], [Gau77a], [Gau77b], [MP21], [Pop13] and the references therein. Arguably, the two most heavily studied ones are the {\bf SKT} and the {\bf balanced} metrics/manifolds. (We refer to the next section for the relevant definitions.) 

In this paper, we study these two classes of complex manifolds and, partially inspired by [DP20], we take up the question of characterising K\"ahlerianity within them. Indeed, it is well known that there are many compact complex manifolds which admit SKT but no K\"ahler metrics. For example, any compact complex surface is SKT, whereas K\"ahlerianity of a surface is a topological property, equivalent to the first Betti number being even. Meanwhile, starting from dimension 3, there are numerous examples of non-K\"ahler balanced manifolds. 

Nevertheless, the question of whether a compact complex manifold that is both SKT and balanced must be K\"ahler (often referred to in the literature as the {\it Fino-Vezzoni conjecture}, see [FV15]) is still open. It is currently only known to have an affirmative answer in some special cases, such as for a certain class of complex structures on connected sums of several copies of $S^3\times S^3$ ([FLY12]) and on compact complex nilmanifolds endowed with invariant complex structures ([AN22], [FV16]). Specifically, Arroyo and Nicolini proved in [AN22] that if such a nilmanifold admits an SKT metric, then it must be $k$-step nilpotent with $k\leq 2$. This implies the validity of the Fino-Vezzoni conjecture on these manifolds.\footnote{The authors are very grateful to L. Ugarte for bringing the reference [AN22] and the results therein to their attention.} Meanwhile, the conjecture has implicitly been proved for twistor spaces (known to always carry balanced metrics) by Verbitsky in [Ver14] and for {\it class ${\cal C}$ manifolds} by Chiose in [Chi14].

Independently of the above-mentioned classes of manifolds, this conjecture was given an affirmative answer very recently in [Pop22] for general manifolds satisfying extra assumptions, including the {\it pluriclosed star split} hypothesis (a notion introduced therein) on a pair of Hermitian metrics. We hope that one of the future applications of our results of the present paper will be a further line of attack on the general case of the Fino-Vezzoni conjecture as pointed out in $\S$\ref{section:simultaneous_SKT-bal}. 
 
Our approach is variational: we propose two new functionals defined on the spaces of SKT, respectively balanced, metrics and prove that their critical points are the K\"ahler metrics, if any. As these functionals are non-negative by definition, a family of SKT metrics (respectively a family of balanced metrics) asymptotically minimising the corresponding functional {\bf either} converges to a K\"ahler metric {\bf or} degenerates. After a suitable normalisation of the metrics, this should produce in the limit an SKT {\bf $(1,1)$-current}, resp. a balanced {\bf $(n-1,n-1)$-current}. Such a current ought to constitute an {\bf obstruction} to the existence of K\"ahler metrics (see discussion in $\S$\ref{subsection:normalised_SKT-functionals}). We expect that a suitable pairing of two such currents, one SKT, the other balanced, on an SKT balanced manifold will help with the efforts of settling the K\"ahlerianity conjecture for the SKT balanced manifolds.

The main idea behind the functionals is as follows: in the SKT case the pluriclosedness of the metric $\omega$ can be rewritten as $\partial\omega\in ker d$ i.e. $\partial\omega$ is {\bf orthogonal} to $Im\ d^\star$. Exploiting the three-space decomposition at the level of complex-valued $3$-forms, this means that $\partial\omega$ can be decomposed into a harmonic part and a $d$-exact part. The harmonic subspace being always finite-dimensional (by standard elliptic theory on a compact manifold), we consider the orthogonal projection onto the image of $d$ and, generalising a construction from [DP20], we define a {\bf torsion form} $\rho_\omega$ as the minimal $L^2$-norm primitive (for the operator $d$) of the projection of $\partial\omega$. The functional $F$ is then defined as the squared $L^2$-norm of this torsion form. The major difference from [DP20], where we dealt with Hermitian-symplectic metrics, is that $\rho_\omega$ is now merely a (complex-valued) $2$-form and its decomposition into pure types plays a major role in the computations. 

Similarly, in the balanced case we exploit the three-space decomposition associated with the $\bar\partial$-Laplacian on $(n-1,n-1)$-forms and the orthogonal projection onto $Im\ \bar\partial$ and again we define a suitable torsion form which is now of pure type $(n-1,n-2)$.

In this paper, we mainly deal with non-normalised functionals, but we also include a discussion of normalised versions of them in $\S$\ref{subsection:normalised_SKT-functionals}. Our main results consist in the computations of the critical point equations and the variation of our functionals in the scaling direction. To make this possible, we first compute the first variation of various differential geometric operators that are frequently used in this context, such as the Hodge star operator $\star_\omega$, the adjoint $\Lambda_\omega$ of the Lefschetz operator $L_\omega$ of multiplication by a metric $\omega$ and the orthogonal projectors onto the kernels of the Laplacians $\Delta'$, $\Delta''$, $\Delta$. We hope that these preliminary computations are of sufficient generality to be of independent interest. 

As our functionals are defined on open cones, it is natural to expect that if minimisers thereof exist, the variation in the scaling direction has to vanish and the normalised functional has to be constant along the rays. This is confirmed by our computations which, moreover, show that vanishing of any of the functionals yields a K\"ahler metric.

The paper is organised as follows: in Section \ref{section:two-functionals}, after fixing the notation and recalling a few basic definitions, we define the functional $F$ acting on the space of SKT metrics on a fixed compact complex manifold and then the functional $G$ acting on $(n-1)$-st wedge powers of balanced metrics. Section \ref{section:1st-variation_operators} is devoted to the computations of the first variations of various operators that are involved in the definitions of $F$ and $G$. The next two sections contain the computations of the critical points of both functionals (Section \ref{section:critical-points-computation_functional_F} for $F$ and Section \ref{section:critical-points-computation_functional_G} for $G$) and the variations in the scaling direction. We utilise here the computations of Section \ref{section:1st-variation_operators} together with a trick that we already used in [DP20]. In the last section \ref{section:simultaneous_SKT-bal}, we briefly discuss a simultaneous approach to the variations of SKT and balanced metrics, when both are supposed to exist, on a compact complex manifold by proposing a third functional $H$ depending on two arguments. The starting point is an observation that generalises the quick proof given in [Pop15] to the fact that if a same Hermitian metric is both SKT and balanced then it must be K\"ahler. In $\S$\ref{subsection:normalised_SKT-functionals}, we discuss certain normalisations of our SKT functional $F$ and the role they may play in constructing new obstructions to compact SKT manifolds admitting K\"ahler metrics. 

\vspace{2ex}

\noindent {\it Acknowledgements}.  The first-named author is partially supported by grant no. 2021/41/B/ST1/01632 from the National Science Center, Poland. Both authors are very grateful to the referee and to L. Ugarte for useful comments and for pointing out references the authors were previously unaware of. 

\section{Two functionals}\label{section:two-functionals}  

Throughout the paper $X$ will be a fixed compact complex manifold with $\mbox{dim}_\C X=n$. We tacitly assume $n\geq 2$, although some of the computations also make sense for compact Riemann surfaces. As is customary, we identify Hermitian metrics on $X$ with the corresponding $C^\infty$ positive definite $(1,1)$-forms. For any $(1,\,1)$-form $\gamma$ on $X$ and any $p\in\{1,\dots , n\}$, we use the notation $\gamma_p:=\gamma^p/p!$. 

Recall that a metric $\omega$ is said to be {\bf SKT} if $\partial\bar\partial\omega=0$. A manifold admitting such a metric is called an SKT manifold. These metrics are also known as {\it pluriclosed} in the literature (see [ST10]). A metric $\omega$ is said to be {\bf balanced} if $d\omega^{n-1}=0$ (see [Gau77b] for the notion, [Mic83] for the name) and again any manifold admitting such metrics is called a balanced manifold. Both notions impose topological and analytical restrictions on the underlying manifold -- see [Pop15] for more details. It is by now a standard fact (see [AI01] or [IP12] for a proof or [Pop15] for a short proof) that a metric which is both SKT and balanced is K\"ahler. 

We will study conditions under which  these two special types of Hermitian metrics have the stronger K\"ahler property. We introduce a functional in each of these two cases and prove that its critical points are precisely the K\"ahler metrics on $X$.

\subsection{The SKT functional}\label{subsection:functional_skt-Delta}

Suppose there exists an {\bf SKT metric} $\omega$ on $X$, namely a $C^\infty$ positive definite $(1,\,1)$-form $\omega$ on $X$ such that $\partial\bar\partial\omega=0$. The last condition is equivalent to $\partial\omega\in\ker d$, hence also to $$\star(\partial\omega)\in\ker d^\star,$$ where $\star=\star_\omega$ is the Hodge star operator induced by $\omega$ and $d^\star = d^\star_\omega$ is the adjoint of $d$ w.r.t. the $L^2$-inner product $\langle\langle\,\,,\,\,\rangle\rangle_\omega$ defined by $\omega$. We used the facts that $d^\star = -\star d\star$, $\star$ is an isomorphism and $\star\star = -\mbox{Id}$ on odd-degreed forms. (Note that $\partial\omega$ is of degree $3$.)

Now, let $P^{(3)}_\omega:C^\infty_3(X,\,\C)\longrightarrow\mbox{Im}\,d$ be the orthogonal projection (w.r.t. the $L^2$-inner product defined by $\omega$) onto the subspace $\mbox{Im}\,d$ of smooth $3$-forms, induced by the classical orthogonal $3$-space decomposition: \begin{equation}\label{eqn:3-space-decomp_d-3}C^\infty_3(X,\,\C) = {\cal H}^3_\Delta(X,\,\C)\oplus\mbox{Im}\,d\oplus\mbox{Im}\,d^\star\end{equation} in which ${\cal H}^3_\Delta(X,\,\C)$ is the kernel of the $d$-Laplacian $\Delta=\Delta_\omega=dd^\star_\omega + d^\star_\omega d:C^\infty_3(X,\,\C)\longrightarrow C^\infty_3(X,\,\C)$ induced by $\omega$, $\ker d = {\cal H}^3_\Delta(X,\,\C)\oplus\mbox{Im}\,d$ and $\ker d^\star = {\cal H}^3_\Delta(X,\,\C)\oplus\mbox{Im}\,d^\star$.

Since $\partial\omega\in\ker d$, the projection of $\partial\omega\in C^\infty_3(X,\,\C)$ onto $ \mbox{Im}\,d^\star$ vanishes.

\begin{Def}\label{Def:torsion-form} Given an {\bf SKT} metric $\omega$ on a compact complex manifold $X$, the unique smooth $2$-form $\rho=\rho_\omega\in C^\infty_2(X,\,\C)$ that solves the system of equations: \begin{eqnarray}\label{eqn:system_torsion-form}(a)\hspace{1ex} d\rho_\omega = P^{(3)}_\omega(\partial\omega) \hspace{3ex} \mbox{and} \hspace{3ex} (b)\hspace{1ex}\rho_\omega\in\mbox{Im}\,d^\star_\omega\end{eqnarray} is called the {\bf torsion form} of $\omega$.

\end{Def}  

Note that equation (a) of (\ref{eqn:system_torsion-form}) is solvable since $P^{(3)}_\omega(\partial\omega)\in\mbox{Im}\,d$. However, the solution of equation (a) is unique only up to $\ker d$, so the extra condition (b) makes it unique in the absolute sense since $\mbox{Im}\,d^\star$ is the orthogonal complement w.r.t. $\langle\langle\,\,,\,\,\rangle\rangle_\omega$ of $\ker d$ inside the space $C^\infty_2(X,\,\C)$. For the same reason, the solution $\rho_\omega$ of system (\ref{eqn:system_torsion-form}) is the {\it minimal $L^2_\omega$-norm solution} of equation (a). It is given by the classical {\it Neumann formula}: \begin{eqnarray}\label{eqn:Neumann_system_torsion-form}\rho_\omega = \Delta_\omega^{-1}d^\star_\omega(P^{(3)}_\omega(\partial\omega)),\end{eqnarray} where $\Delta_\omega^{-1}$ is the {\it Green operator} of the (elliptic) $d$-Laplacian $\Delta_\omega=dd^\star_\omega + d^\star_\omega d:C^\infty_2(X,\,\C)\longrightarrow C^\infty_2(X,\,\C)$ induced by $\omega$.

\begin{Lem}\label{Lem:Kaehler_torsion-vanishing} Let $\omega$ be an {\bf SKT} metric on a compact complex manifold $X$. The following equivalences hold: \begin{eqnarray*}\label{eqn:Kaehler_torsion-vanishing}(i)\hspace{1ex} \omega \hspace{1ex}\mbox{is K\"ahler}\hspace{1ex} \iff \hspace{1ex} (ii)\hspace{1ex}\rho_\omega =0 \hspace{1ex} \iff \hspace{1ex} (iii) \hspace{1ex} P^{(3)}_\omega(\partial\omega)=0.\end{eqnarray*}

\end{Lem}  

\noindent {\it Proof.} We will prove the implications: $(i)\implies (ii) \implies (iii) \implies (i)$.

\hspace{1ex}

$(i)\implies (ii).$ If $\omega$ is K\"ahler, $\partial\omega=0$, hence $P^{(3)}_\omega(\partial\omega)=0$, so the minimal $L^2$-norm solution of equation $d\rho_\omega = 0$ is $\rho_\omega=0$. 

\hspace{1ex}

$(ii)\implies (iii).$ If $\rho_\omega=0$, then $P^{(3)}_\omega(\partial\omega)=d\rho_\omega=0$.

\hspace{1ex}

$(iii)\implies (i).$ Suppose that $P^{(3)}_\omega(\partial\omega)=0$. Since $\partial\omega\in\ker d = {\cal H}_\Delta(X,\,\C)\oplus\mbox{Im}\,d$, this is equivalent to $\partial\omega\in{\cal H}_\Delta(X,\,\C)$.

On the other hand, for any pure-type form $\alpha$ we have: \begin{eqnarray}\label{eqn:Laplacians_pure-type}\langle\langle\Delta\alpha,\,\alpha\rangle\rangle = \langle\langle\Delta'\alpha,\,\alpha\rangle\rangle + \langle\langle\Delta''\alpha,\,\alpha\rangle\rangle\end{eqnarray} Indeed, we have: \begin{eqnarray*}\langle\langle\Delta\alpha,\,\alpha\rangle\rangle = ||d\alpha||^2 + ||d^\star\alpha||^2 \stackrel{(a)}{=} ||\partial\alpha||^2 + ||\partial^\star\alpha||^2 + ||\bar\partial\alpha||^2 + ||\bar\partial^\star\alpha||^2 = \langle\langle\Delta'\alpha,\,\alpha\rangle\rangle + \langle\langle\Delta''\alpha,\,\alpha\rangle\rangle,\end{eqnarray*} where the crucial equality (a) follows from $\partial\alpha$ and $\bar\partial\alpha$ being of {\it different} pure types, hence orthogonal to each other, hence $||d\alpha||^2 = ||\partial\alpha + \bar\partial\alpha||^2 = ||\partial\alpha||^2 + ||\bar\partial\alpha||^2$, while $\partial^\star\alpha$ and $\bar\partial^\star\alpha$ too are of {\it different} pure types, hence orthogonal to each other, hence $||d^\star\alpha||^2 = ||\partial^\star\alpha + \bar\partial^\star\alpha||^2 = ||\partial^\star\alpha||^2 + ||\bar\partial^\star\alpha||^2$.

  Now, taking $\alpha=\partial\omega$ (which is of pure type $(2,\,1)$) in (\ref{eqn:Laplacians_pure-type}), the property $\Delta(\partial\omega)=0$ implies that $\Delta'(\partial\omega)=0$ and $\Delta''(\partial\omega)=0$. In particular, we get $$\partial\omega\in\ker\Delta'\cap\mbox{Im}\,\partial =\{0\},$$ where the last equality of spaces follows from $\ker\Delta'\perp\mbox{Im}\,\partial$.

Thus, $\partial\omega=0$, which means that $\omega$ is K\"ahler.  \hfill $\Box$

\vspace{2ex}

Now, we consider the set \begin{equation*}\label{eqn:def-Somega}{\cal S}:=\bigg\{\omega\,\mid\,\omega\in C^\infty_{1,\,1}(X,\,\R) \hspace{1ex}\mbox{such that}\hspace{1ex} \omega>0 \hspace{1ex}\mbox{and}\hspace{1ex} \partial\bar\partial\omega= 0\bigg\}\subset \ker(\partial\bar\partial)\cap C^\infty_{1,\,1}(X,\,\R)\end{equation*} of all the SKT metrics on $X$. The set ${\cal S}$ is an {\it open convex cone} in the real vector space $\ker(\partial\bar\partial)\cap C^\infty_{1,\,1}(X,\,\R)$.

\begin{Def}\label{Def:F_energy-functional_SKT} Let $X$ be a compact complex SKT manifold with $\mbox{dim}_{\C}X=n$. We define the following {\bf energy functional}: \begin{equation}\label{eqn:F_energy-functional_SKT}F : {\cal S} \to [0,\,+\infty), \hspace{3ex} F(\omega) = \int\limits_X|\rho_\omega|^2_\omega\,dV_\omega = ||\rho_\omega||^2_\omega,  \end{equation}

\noindent where $\rho_{\omega}$ is the torsion form of the SKT metric $\omega\in{\cal S}$ defined in Definition \ref{Def:torsion-form}, while $|\,\,\,|_\omega$ is the pointwise norm and $||\,\,\,||_\omega$ is the $L^2$-norm induced by $\omega$.

\end{Def}

The first trivial observation that justifies the introduction of the functional $F$ is the following.

\begin{Lem}\label{Lem:vanishing-F-Kaehler} Let $\omega$ be an SKT metric on $X$. The following equivalence holds:

\begin{equation}\label{eqn:vanishing-F-Kaehler}\omega \hspace{1ex} \mbox{is K\"ahler} \iff F(\omega)=0.\end{equation}

\end{Lem} 

\noindent {\it Proof.} The condition $F(\omega)=0$ is equivalent to $\rho_{\omega}=0$. This, in turn, is equivalent to $\omega$ being K\"ahler by Lemma \ref{Lem:Kaehler_torsion-vanishing}. \hfill $\Box$

\vspace{3ex}

The next observation goes further after specifying the behaviour of $F(\omega)$ under rescalings of $\omega$. 

\begin{Lem}\label{Lem:rescaling_F-critical} Let $\omega$ be an SKT metric on $X$. For every real $\lambda>0$, the following equality holds: \begin{equation}\label{eqn:rescaling_F-critical}F(\lambda\,\omega) = \lambda^n\,F(\omega).\end{equation}

In particular, \begin{eqnarray}\label{eqn:rescaling_F-critical_1+t}\frac{d}{dt}_{|t=0}\bigg(F((1+t)\omega)\bigg) = nF(\omega)\end{eqnarray} and the following equivalence holds:

\begin{equation}\label{eqn:critical-points_F-Kaehler}\omega \hspace{1ex} \mbox{is K\"ahler} \iff \omega \hspace{1ex} \mbox{is a critical point for} \hspace{1ex} F.\end{equation}

\end{Lem} 

\noindent {\it Proof.} The equivalence (\ref{eqn:critical-points_F-Kaehler}) follows from Lemma \ref{Lem:vanishing-F-Kaehler} and from (\ref{eqn:rescaling_F-critical_1+t}) by considering the vanishing of $(d/dt)_{|t=0}F(\omega+t\gamma)$ with the tangent direction provided by the choice $\gamma=\omega$. So, only (\ref{eqn:rescaling_F-critical}) needs a proof.

Lemma 3.6. in [BP18] shows that, for any Hermitian metric $\omega$ and any positive real $\lambda$, one has $\bar\partial^\star_{\lambda\omega} = (1/\lambda)\,\bar\partial^\star_\omega$ and $\Delta''_{\lambda\omega} = (1/\lambda)\,\Delta''_\omega$. The same argument yields: \begin{eqnarray}\label{eqn:rescaling_d-star_Delta}d^\star_{\lambda\omega} = \frac{1}{\lambda}\,d^\star_\omega \hspace{3ex} \mbox{and} \hspace{3ex} \Delta_{\lambda\omega} = \frac{1}{\lambda}\,\Delta_\omega\end{eqnarray} on differential forms of any degree. In particular, $\Delta_\omega$ and $\Delta_{\lambda\omega}$ have the same kernels, while $d^\star_{\lambda\omega}$ and $d^\star_\omega$ have the same image.

Now, let $P_{\Delta,\,\omega}$, resp. $P_{\Delta,\,\lambda\,\omega}$, be the $L^2_\omega$-orthogonal projection, resp. the $L^2_{\lambda\omega}$-orthogonal projection, from the space of $C^\infty$ $3$-forms onto the kernel of $\Delta_\omega$, resp. the kernel of $\Delta_{\lambda\omega}$. The decomposition $\partial\omega = P_{\Delta,\,\omega}(\partial\omega) + d\rho_\omega$, with $\rho_\omega\in\mbox{Im}\,d^\star_\omega$, leads to $\partial(\lambda\,\omega) = \lambda\,P_{\Delta,\,\omega}(\partial\omega) + d(\lambda\rho_\omega)$ and, thanks to (\ref{eqn:rescaling_d-star_Delta}), the latter is the decomposition of $\partial(\lambda\,\omega)\in\ker d$ according to the $L^2_{\lambda\omega}$-orthogonal splitting $\ker d = \ker\Delta_{\lambda\omega}\oplus\mbox{Im}\,d$. Moreover, $\lambda\rho_\omega\in\mbox{Im}\,d^\star_\omega = \mbox{Im}\,d^\star_{\lambda\omega}$, so we get \begin{eqnarray*}\rho_{\lambda\omega} = \lambda\rho_\omega.\end{eqnarray*}

From this, we infer: \begin{eqnarray*}F(\lambda\omega) = \int\limits_X|\rho_{\lambda\omega}|^2_{\lambda\omega}\,dV_{\lambda\omega} = \lambda^2\int\limits_X\frac{1}{\lambda^2}\,|\rho_\omega|^2_\omega\,(\lambda^n\,dV_\omega) = \lambda^nF(\omega), \hspace{5ex} \lambda>0,\end{eqnarray*} which proves (\ref{eqn:rescaling_F-critical}).  \hfill $\Box$

\vspace{2ex}

In view of further use, we note that the above proof also gives $P^{(3)}_{\lambda\omega}(\partial(\lambda\omega)) = \lambda\,P^{(3)}_\omega(\partial\omega)$, hence \begin{eqnarray}\label{eqn:P^3_del-omega_lambda}P^{(3)}_{\lambda\omega}(\partial\omega) = P^{(3)}_\omega(\partial\omega),  \hspace{5ex} \lambda>0.\end{eqnarray}

\subsection{The balanced functional}\label{subsection:functional_balanced}

Suppose there exists a balanced metric $\omega$ on $X$, namely a $C^\infty$ positive definite $(1,\,1)$-form $\omega$ such that $\bar\partial\omega_{n-1}=0$. (Recall the notation $\omega_{n-1}:=\omega^{n-1}/(n-1)!$, where $n=\mbox{dim}_\C X$.) Since every $C^\infty$ positive definite $(n-1,\,n-1)$-form $\Omega$ on $X$ admits a unique $(n-1)$-st root, namely a unique $C^\infty$ positive definite $(1,\,1)$-form $\omega$ such that $\omega_{n-1} = \Omega$ (see [Mic83]), $\omega$ and $\omega_{n-1}$ determine each other. Thus, we can refer to either $\omega$ or $\omega_{n-1}$ as a balanced metric when $\omega_{n-1}$ lies in $\ker\bar\partial$.

For any balanced metric $\omega$ on $X$, let  $\Delta''=\Delta''_\omega=\bar\partial\bar\partial^\star_\omega + \bar\partial^\star_\omega\bar\partial:C^\infty_{n-1,\,n-1}(X,\,\C)\longrightarrow C^\infty_{n-1,\,n-1}(X,\,\C)$ be the induced $\bar\partial$-Laplacian acting in bidegree $(n-1,\,n-1)$ and let \begin{equation*}C^\infty_{n-1,\,n-1}(X,\,\C) = \ker\Delta_\omega''\oplus\mbox{Im}\,\bar\partial\oplus\mbox{Im}\,\bar\partial^\star\end{equation*} be the classical $L^2_\omega$-orthogonal $3$-space decomposition. We denote by $P^{(n-1,\,n-1)}_{\omega_{n-1}}:C^\infty_{n-1,\,n-1}(X,\,\C)\longrightarrow\mbox{Im}\,\bar\partial$ the $L^2_\omega$-orthogonal projection onto the image of $\bar\partial$. It is standard (for any Hermitian metric $\omega$) that $\ker\bar\partial = \ker\Delta_\omega''\oplus\mbox{Im}\,\bar\partial$, so if $\omega$ is balanced, $\omega_{n-1}\in\ker\bar\partial$ and therefore $\omega_{n-1}$ splits uniquely as an $L^2_\omega$-orthogonal sum of a $\Delta_\omega''$-harmonic form and $P^{(n-1,\,n-1)}_{\omega_{n-1}}(\omega_{n-1})$.

\begin{Def}\label{Def:balanced_torsion-form} Given a {\bf balanced} metric $\omega$ on a compact complex $n$-dimensional manifold $X$, the unique smooth $(n-1,\,n-2)$-form $\Gamma=\Gamma_{\omega_{n-1}}\in C^\infty_{n-1,\,n-2}(X,\,\C)$ that solves the system of equations: \begin{eqnarray}\label{eqn:system_balanced_torsion-form}(a)\hspace{1ex} \bar\partial\Gamma_{\omega_{n-1}} = P^{(n-1,\,n-1)}_{\omega_{n-1}}(\omega_{n-1}) \hspace{3ex} \mbox{and} \hspace{3ex} (b)\hspace{1ex}\Gamma_{\omega_{n-1}}\in\mbox{Im}\,\bar\partial^\star_\omega\end{eqnarray} is called the {\bf torsion $(n-1,\,n-2)$-form} of $\omega_{n-1}$.

\end{Def}

Note that equation (a) of (\ref{eqn:system_balanced_torsion-form}) is solvable since $P^{(n-1,\,n-1)}_{\omega_{n-1}}(\omega_{n-1})\in\mbox{Im}\,\bar\partial$. However, the solution of equation (a) is unique only up to $\ker\bar\partial$, so the extra condition (b) makes it unique in the absolute sense since $\mbox{Im}\,\bar\partial^\star_\omega$ is the orthogonal complement of $\ker\bar\partial$ w.r.t. the $L^2_\omega$-inner product $\langle\langle\,\,,\,\,\rangle\rangle_\omega$ inside the space $C^\infty_{n-1,\,n-1}(X,\,\C)$. For the same reason, the solution $\Gamma_{\omega_{n-1}}$ of system (\ref{eqn:system_balanced_torsion-form}) is the {\it minimal $L^2_\omega$-norm solution} of equation (a). It is given by the classical {\it Neumann formula}: \begin{eqnarray}\label{eqn:Neumann_system_balanced_torsion-form}\Gamma_{\omega_{n-1}} = \Delta_\omega^{''-1}\bar\partial^\star_\omega(P^{(n-1,\,n-1)}_{\omega_{n-1}}(\omega_{n-1})),\end{eqnarray} where $\Delta_\omega^{''-1}$ is the {\it Green operator} of the (elliptic) $\bar\partial$-Laplacian $\Delta_\omega''$.

\begin{Lem}\label{Lem:Kaehler_balanced_torsion-vanishing} Let $\omega$ be a {\bf balanced} metric on a compact complex $n$-dimensional manifold $X$. The following equivalences hold: \begin{eqnarray*}\label{eqn:Kaehler_balanced_torsion-vanishing}(i)\hspace{1ex} \omega \hspace{1ex}\mbox{is K\"ahler}\hspace{1ex} \iff \hspace{1ex} (ii)\hspace{1ex} P^{(n-1,\,n-1)}_{\omega_{n-1}}(\omega_{n-1})=0 \hspace{1ex} \iff \hspace{1ex} (iii) \hspace{1ex} \Gamma_{\omega_{n-1}}=0.\end{eqnarray*}

\end{Lem}  

\noindent {\it Proof.} We will prove the implications: $(i)\implies (ii) \implies (iii) \implies (i)$.

\hspace{1ex}

$(i)\implies (ii).$ If $\omega$ is K\"ahler, $\omega_{n-1}\in\ker\Delta_\omega''$, hence $P^{(n-1,\,n-1)}_{\omega_{n-1}}(\omega_{n-1})=0$. Indeed, since $\ker\Delta_\omega'' = \ker\bar\partial\cap\ker\bar\partial^\star$, to see that $\omega_{n-1}$ lies in $\ker\Delta_\omega''$, we notice that:

\vspace{1ex}

$\bar\partial\omega_{n-1} = \omega_{n-2}\wedge\bar\partial\omega = 0$ since $\bar\partial\omega = 0$ by the K\"ahler assumption on $\omega$;

\vspace{1ex}

$\bar\partial^\star\omega_{n-1} = -\star\partial(\star\omega_{n-1}) = -\star\partial\omega = 0$ since $\partial\omega = 0$ by the K\"ahler assumption on $\omega$.

\hspace{1ex}

$(ii)\implies (iii).$ If $P^{(n-1,\,n-1)}_{\omega_{n-1}}(\omega_{n-1})=0$, the zero form is a solution of equation (a) of (\ref{eqn:system_balanced_torsion-form}). Then, the zero form is necessarily the minimal $L^2_\omega$-norm solution $\Gamma_{\omega_{n-1}}$ of that equation, hence $\Gamma_{\omega_{n-1}}=0$.

\hspace{1ex}

$(iii)\implies (i).$ Suppose that $\Gamma_{\omega_{n-1}}=0$. Then $\bar\partial\Gamma_{\omega_{n-1}}=0$, hence $P^{(n-1,\,n-1)}_{\omega_{n-1}}(\omega_{n-1})=0$. Since $\omega_{n-1}\in\ker\bar\partial = \ker\Delta_\omega''\oplus\mbox{Im}\,\bar\partial$, by the balanced hypothesis on $\omega$, the property $P^{(n-1,\,n-1)}_{\omega_{n-1}}(\omega_{n-1})=0$ is equivalent to $\omega_{n-1}\in\ker\Delta_\omega'' = \ker\bar\partial\cap\ker\bar\partial^\star$. In particular, $\bar\partial^\star\omega_{n-1}=0$, which means that $0=-\star\partial(\star\omega_{n-1}) = -\star\partial\omega$. Since $\star$ is an isomorphism, this amounts to $\partial\omega=0$, which means that $\omega$ is K\"ahler.  \hfill $\Box$

\vspace{2ex}

Now, we consider the set \begin{equation*}\label{eqn:def-Bomega}{\cal B}:=\bigg\{\omega_{n-1}\,\mid\,\omega\in C^\infty_{1,\,1}(X,\,\R) \hspace{1ex}\mbox{such that}\hspace{1ex} \omega>0 \hspace{1ex}\mbox{and}\hspace{1ex} \bar\partial\omega_{n-1}= 0\bigg\}\subset \ker(\bar\partial)\cap C^\infty_{n-1,\,n-1}(X,\,\R)\end{equation*} of all the balanced metrics on $X$. The set ${\cal B}$ is an {\it open convex cone} in the real vector space $\ker(\bar\partial)\cap C^\infty_{n-1,\,n-1}(X,\,\R)$.

\begin{Def}\label{Def:F_energy-functional_balanced} Let $X$ be a compact complex balanced manifold with $\mbox{dim}_{\C}X=n$. We define the following {\bf energy functional}: \begin{equation}\label{eqn:G_energy-functional_balanced}G :{\cal B} \to [0,\,+\infty), \hspace{3ex} G(\omega_{n-1}) = \int\limits_X|\Gamma_{\omega_{n-1}}|^2_\omega\,dV_\omega = ||\Gamma_{\omega_{n-1}}||^2_\omega,  \end{equation}

\noindent where $\Gamma_{\omega_{n-1}}$ is the torsion $(n-1,\,n-2)$-form of the balanced metric $\omega_{n-1}\in{\cal B}$ defined in Definition \ref{Def:balanced_torsion-form}, while $|\,\,\,|_\omega$ is the pointwise norm and $||\,\,\,||_\omega$ is the $L^2$-norm induced by $\omega$.

\end{Def}

The first trivial observation that justifies the introduction of the functional $G$ is the following.

\begin{Lem}\label{Lem:vanishing-G-Kaehler} Let $\omega$ be a balanced metric on $X$. The following equivalence holds:

\begin{equation}\label{eqn:vanishing-G-Kaehler}\omega \hspace{1ex} \mbox{is K\"ahler} \iff G(\omega_{n-1})=0.\end{equation}

\end{Lem} 

\noindent {\it Proof.} The condition $G(\omega_{n-1})=0$ is equivalent to $\Gamma_{\omega_{n-1}}=0$. This, in turn, is equivalent to $\omega$ being K\"ahler by Lemma \ref{Lem:Kaehler_balanced_torsion-vanishing}. \hfill $\Box$

\vspace{3ex}

The next observation is the analogue in the balanced case of Lemma \ref{Lem:rescaling_F-critical}.

\begin{Lem}\label{Lem:rescaling_G-critical} Let $\omega$ be a balanced metric on $X$. For every real $\lambda>0$, the following equality holds: \begin{equation}\label{eqn:rescaling_G-critical}G(\lambda\,\omega_{n-1}) = \lambda^{\frac{n+1}{n-1}}\,G(\omega_{n-1}).\end{equation}

In particular, \begin{eqnarray}\label{eqn:rescaling_G-critical_1+t}\frac{d}{dt}_{|t=0}\bigg(G((1+t)\omega_{n-1})\bigg) = \frac{n+1}{n-1}\,G(\omega_{n-1})\end{eqnarray} and the following equivalence holds:

\begin{equation}\label{eqn:critical-points_G-Kaehler}\omega \hspace{1ex} \mbox{is K\"ahler} \iff \omega_{n-1} \hspace{1ex} \mbox{is a critical point for} \hspace{1ex} G.\end{equation}

\end{Lem} 

\noindent {\it Proof.} The equivalence (\ref{eqn:critical-points_G-Kaehler}) follows from Lemma \ref{Lem:vanishing-G-Kaehler} and from (\ref{eqn:rescaling_G-critical_1+t}) by considering the vanishing of $(d/dt)_{|t=0}G(\omega_{n-1}+t\Omega)$ with the tangent direction provided by the choice $\Omega=\omega_{n-1}$. So, only (\ref{eqn:rescaling_G-critical}) needs a proof.

As previously said, Lemma 3.6. in [BP18] shows that, for any Hermitian metric $\omega$ and any positive real $\lambda$, one has: \begin{eqnarray}\label{eqn:rescaling_d-bar-star_Delta''}\bar\partial^\star_{\lambda\omega} = \frac{1}{\lambda}\,\bar\partial^\star_\omega \hspace{3ex} \mbox{and} \hspace{3ex} \Delta''_{\lambda\omega} = \frac{1}{\lambda}\,\Delta''_\omega\end{eqnarray} on differential forms of any bidegree. In particular, $\Delta''_\omega$ and $\Delta''_{\lambda\omega}$ have the same kernels, while $\bar\partial^\star_{\lambda\omega}$ and $\bar\partial^\star_\omega$ have the same image.

Now, let $P_{\Delta'',\,\omega}$, resp. $P_{\Delta'',\,\lambda\,\omega}$, be the $L^2_\omega$-orthogonal projection, resp. the $L^2_{\lambda\omega}$-orthogonal projection, from the space of $C^\infty$ $(n-1,\,n-1)$-forms onto the kernel of $\Delta''_\omega$, resp. the kernel of $\Delta''_{\lambda\omega}$. The decomposition $\omega_{n-1} = P_{\Delta'',\,\omega}(\omega_{n-1}) + \bar\partial\Gamma_{\omega_{n-1}}$, with $\Gamma_{\omega_{n-1}}\in\mbox{Im}\,\bar\partial^\star_\omega$, leads to $(\lambda^{\frac{1}{n-1}}\,\omega)_{n-1} = \lambda\,P_{\Delta'',\,\omega}(\omega_{n-1}) + \bar\partial(\lambda\Gamma_{\omega_{n-1}})$ and, thanks to (\ref{eqn:rescaling_d-bar-star_Delta''}), the latter is the decomposition of $\lambda\,\omega_{n-1}\in\ker\bar\partial$ according to the $L^2_{\lambda\omega}$-orthogonal splitting $\ker\bar\partial = \ker\Delta''_{\lambda\omega}\oplus\mbox{Im}\,\bar\partial$. Moreover, $\lambda\Gamma_{\omega_{n-1}}\in\mbox{Im}\,\bar\partial^\star_\omega = \mbox{Im}\,\bar\partial^\star_{\lambda^{1/(n-1)}\omega}$, so we get \begin{eqnarray*}\Gamma_{\lambda\omega_{n-1}} = \lambda\Gamma_{\omega_{n-1}}.\end{eqnarray*}

From this, we infer: \begin{eqnarray*}G(\lambda\omega_{n-1}) = \int\limits_X|\Gamma_{\lambda\omega_{n-1}}|^2_{\lambda^{1/(n-1)}\omega}\,dV_{\lambda^{1/(n-1)}\omega} & = & \lambda^2\int\limits_X\bigg(\frac{1}{\lambda^{\frac{1}{n-1}}}\bigg)^{2n-3}\,|\Gamma_{\omega_{n-1}}|^2_\omega\,(\lambda^{\frac{n}{n-1}}\,dV_\omega) \\
 & = & \lambda^{\frac{n+1}{n-1}}G(\omega_{n-1}), \hspace{5ex} \lambda>0,\end{eqnarray*} which proves (\ref{eqn:rescaling_G-critical}).  \hfill $\Box$

\vspace{2ex}

In view of further use, we note that the above proof also gives $P^{(n-1,\,n-1)}_{\lambda\omega_{n-1}}(\lambda\omega_{n-1}) = \lambda\,P^{(n-1,\,n-1)}_{\omega_{n-1}}(\omega_{n-1})$, hence \begin{eqnarray}\label{eqn:P^n-1n-1_omega_n-1_lambda}P^{(n-1,\,n-1)}_{\lambda\omega_{n-1}}(\omega_{n-1}) = P^{(n-1,\,n-1)}_{\omega_{n-1}}(\omega_{n-1}),  \hspace{5ex} \lambda>0.\end{eqnarray}

\vspace{2ex}

 We will need to compute the critical points of the energy functionals $F$ and $G$. Before proceeding, we pause to compute the differentials of several of the classical operators involved in our construction that will come in handy later on and have an interest of their own.

 \section{The first variation of various operators }\label{section:1st-variation_operators} Let $X$ be a complex manifold with $\mbox{dim}_\C X=n$. We fix a Hermitian metric $\omega$ on $X$. This means that $\omega$ is a $C^\infty$ positive definite $(1,\,1)$-form on $X$. Meanwhile, let $\gamma$ be a {\it real} (in general signless) $C^\infty$ $(1,\,1)$-form on $X$. We will be considering reals $t$ sufficiently close to $0$ such that $\omega + t\gamma$ is positive definite.

 Fix a point $x\in X$ and choose local holomorphic coordinates $z_1,\dots , z_n$ about $x$ such that \begin{equation}\label{eqn:omega-gamma_diagonalisation}\omega(x)=\sum\limits_{j=1}^n idz_j\wedge d\bar{z}_j  \hspace{3ex} \mbox{and} \hspace{3ex} \gamma(x)=\sum\limits_{j=1}^n \gamma_j\,idz_j\wedge d\bar{z}_j,\end{equation} where $\gamma_1,\dots , \gamma_n\in\R$ are the eigenvalues of $\gamma$ w.r.t. $\omega$ at $x$. Then, obviously, \begin{equation}\label{eqn:omega-plus-t-gamma_x}(\omega + t\gamma)(x)=\sum\limits_{j=1}^n (1+ t\gamma_j)\, idz_j\wedge d\bar{z}_j,  \hspace{6ex} t\in\R.\end{equation}

 \subsection{Computation of the first variation of the trace of $(1,\,1)$-forms}\label{subsection:1st-variation_trace_11}

 In the setting described above, let $(\alpha_t)_t$ be a $C^\infty$ family of smooth $(1,\,1)$-forms on $X$. In the local coordinates about a given point $x\in X$ chosen above, we have \begin{equation*}\alpha_t =\sum\limits_{l,\,r=1}^n\alpha^{(t)}_{l\bar{r}}\,idz_l\wedge d\bar{z}_r,\end{equation*} where the $\alpha^{(t)}_{l\bar{r}}$'s are $C^\infty$ functions on a neighbourhood of $x$. For any Hermitian metric $\rho$ on $X$, $\Lambda_\rho$ will denote the adjoint of the multiplication operator $\rho\wedge\cdot$ w.r.t. the pointwise inner product $\langle\,\,,\,\,\rangle_\rho$ defined by $\rho$.

\begin{Lem}\label{Lem:1st-variation_trace_11} For every $C^\infty$ family $(\alpha_t)_{t\in(-\varepsilon,\,\varepsilon)}$ of forms $\alpha_t\in C^\infty_{1,\,1}(X,\,\C)$ with $\varepsilon>0$ so small that $\omega+t\gamma>0$ for all $t\in(-\varepsilon,\,\varepsilon)$, the following formula holds: \begin{eqnarray}\label{eqn:1st-variation_trace_11}\frac{d}{dt}\,(\Lambda_{\omega+t\gamma}\alpha_t) = \Lambda_{\omega+t\gamma}\bigg(\frac{d\alpha_t}{dt} \bigg) -\langle\alpha_t,\,\gamma\rangle_{\omega+t\gamma},  \hspace{10ex} t\in(-\varepsilon,\,\varepsilon).\end{eqnarray}
     
\end{Lem}

 \noindent {\it Proof.} At an arbitrary point $x\in X$ about which the local coordinates have been chosen such that (\ref{eqn:omega-gamma_diagonalisation}) holds, (\ref{eqn:omega-plus-t-gamma_x}) implies the following equality at $x$: \begin{equation*}\Lambda_{\omega+t\gamma}\alpha_t = \sum\limits_{l=1}^n\frac{\alpha^{(t)}_{l\bar{l}}}{1+t\gamma_l},  \hspace{5ex} t\in(-\varepsilon,\,\varepsilon).\end{equation*} Deriving with respect to $t$, we get the first equality below at $x$: \begin{equation*}\frac{d}{dt}\,(\Lambda_{\omega+t\gamma}\alpha_t) = -\sum\limits_{l=1}^n\frac{\alpha^{(t)}_{l\bar{l}}\,\gamma_l}{(1+t\gamma_l)^2} + \sum\limits_{l=1}^n\frac{\frac{d\alpha^{(t)}_{l\bar{l}}}{dt}}{1+t\gamma_l} = -\langle\alpha_t,\,\gamma\rangle_{\omega+t\gamma} + \Lambda_{\omega+t\gamma}\bigg(\frac{d\alpha_t}{dt} \bigg),  \hspace{5ex} t\in(-\varepsilon,\,\varepsilon),\end{equation*} where the second equality follows again from (\ref{eqn:omega-plus-t-gamma_x}).  \hfill $\Box$

 \subsection{Computation of the first variation of the Hodge star operator}\label{subsection:1st-variation_Hodge-star}

 We will prove the following formula for the differential of the Hodge star operator $\star_\omega$ induced by the metric $\omega$ in the direction of the given real $(1,\,1)$-form $\gamma$.

\begin{Lem}\label{Lem:1st-variation_Hodge-star} For every bidegree $(p,\,q)$ and every $C^\infty$ family $(v_t)_{t\in(-\varepsilon,\,\varepsilon)}$ of forms $v_t\in C^\infty_{p,\,q}(X,\,\C)$ with $\varepsilon>0$ so small that $\omega+t\gamma>0$ for all $t\in(-\varepsilon,\,\varepsilon)$, the following formula holds: \begin{eqnarray}\label{eqn:1st-variation_Hodge-star}\frac{d}{dt}\bigg|_{t=0}(\star_{\omega+t\gamma}v_t) = \star_\omega\bigg(\frac{dv_t}{dt}\bigg|_{t=0} + [\Lambda_\omega,\,\gamma\wedge\cdot]\,v_0\bigg),\end{eqnarray} where $\gamma\wedge\cdot$ is the operator of multiplication by $\gamma$ and $\Lambda_\omega$ is the adjoint of the multiplication operator $\omega\wedge\cdot$ w.r.t. the pointwise inner product $\langle\,\,,\,\,\rangle_\omega$ defined by $\omega$.

\end{Lem}

\noindent {\it Proof.} The formula to prove being pointwise, it suffices to prove it at a pregiven point $x$ with a choice of local coordinates such that (\ref{eqn:omega-gamma_diagonalisation}) holds at $x$. 

For every $t\in(-\varepsilon,\,\varepsilon)$, (\ref{eqn:omega-plus-t-gamma_x}) yields the following equalities at $x$: \begin{equation*}\langle dz_j,\,dz_k\rangle_{\omega+t\gamma} = \langle d\bar{z}_j,\,d\bar{z}_k\rangle_{\omega+t\gamma} = \delta_{jk}\,\frac{1}{1+t\gamma_j},  \hspace{6ex} j,k.\end{equation*}

Now, let $u,v_t\in C^\infty_{p,\,q}(X,\,\C)$ be arbitrary forms. In a neighbourhood of $x$, they are of the shape \begin{equation*}u = \sum\limits_{|I|=p \atop |J|=q} u_{I\bar{J}}\,dz_I\wedge d\bar{z}_J  \hspace{3ex} \mbox{and} \hspace{3ex} v_t = \sum\limits_{ |I|=p \atop |J|=q} v_{I\bar{J}}(t)\,dz_I\wedge d\bar{z}_J, \hspace{6ex} t\in(-\varepsilon,\,\varepsilon),\end{equation*} where $I=(i_1<\dots <i_p)$ and $J=(j_1<\dots <j_q)$ are multi-indices of lengths $p$, resp. $q$.

The definition of the Hodge star operator $\omega_{\omega+t\gamma}$ induced by the metric $\omega+t\gamma$ (when the real $t$ is so close to $0$ that $\omega+t\gamma>0$) gives the first equality below at every point, while the second equality holds at $x$ and follows from the above preparations: \begin{eqnarray*}u\wedge\star_{\omega+t\gamma}\bar{v_t} = \langle u,\,v_t\rangle_{\omega+t\gamma}\,\frac{(\omega+t\gamma)^n}{n!} = \sum\limits_{ |I|=p \atop |J|=q} u_{I\bar{J}}\,\overline{v_{I\bar{J}}(t)}\,\bigg(\prod\limits_{l\in I}\frac{1}{1+t\gamma_l}\bigg)\,\bigg(\prod\limits_{r\in J}\frac{1}{1+t\gamma_r} \bigg)\,\frac{(\omega+t\gamma)^n}{n!}.\end{eqnarray*}

Deriving this equality w.r.t. $t$, we get: \begin{eqnarray*}u\wedge\frac{d}{dt}(\star_{\omega+t\gamma}\bar{v_t}) & = & \sum\limits_{ |I|=p \atop |J|=q} u_{I\bar{J}}\,\overline{\frac{dv_{I\bar{J}}(t)}{dt}}\,\bigg(\prod\limits_{l\in I}\frac{1}{1+t\gamma_l}\bigg)\,\bigg(\prod\limits_{r\in J}\frac{1}{1+t\gamma_r} \bigg)\,\frac{(\omega+t\gamma)^n}{n!} \\
  & - & \sum\limits_{ |I|=p \atop |J|=q} u_{I\bar{J}}\,\overline{v_{I\bar{J}}(t)}\bigg(\sum\limits_{l\in I}\frac{\gamma_l}{(1+t\gamma_l)^2}\prod\limits_{s\in I\setminus\{l\}}\frac{1}{1+t\gamma_s}\bigg)\,\bigg(\prod\limits_{r\in J}\frac{1}{1+t\gamma_r} \bigg)\,\frac{(\omega+t\gamma)^n}{n!}  \\
  & - & \sum\limits_{|I|=p \atop |J|=q} u_{I\bar{J}}\,\overline{v_{I\bar{J}}(t)}\,\bigg(\prod\limits_{l\in I}\frac{1}{1+t\gamma_l}\bigg)\,\bigg(\sum\limits_{r\in J}\frac{\gamma_r}{(1+t\gamma_r)^2}\prod\limits_{s\in J\setminus\{r\}}\frac{1}{1+t\gamma_s}\bigg)\,\frac{(\omega+t\gamma)^n}{n!}  \\
  & + & \sum\limits_{ |I|=p \atop |J|=q} u_{I\bar{J}}\,\overline{v_{I\bar{J}}(t)}\,\bigg(\prod\limits_{l\in I}\frac{1}{1+t\gamma_l}\bigg)\,\bigg(\prod\limits_{r\in J}\frac{1}{1+t\gamma_r} \bigg)\,\frac{(\omega+t\gamma)^{n-1}}{(n-1)!}\wedge\gamma.\end{eqnarray*} Taking $t=0$ and using the fact that  \begin{eqnarray*}\frac{\omega^{n-1}}{(n-1)!}\wedge\gamma = (\Lambda_\omega\gamma)\,\frac{\omega^n}{n!} = (\gamma_1+\dots + \gamma_n)\,\frac{\omega^n}{n!},\end{eqnarray*} with the last equality holding at $x$ only, we get the first equality below at $x$: \begin{eqnarray*}u\wedge\frac{d}{dt}\bigg|_{t=0}(\star_{\omega+t\gamma}\bar{v_t}) & = & \bigg\langle u,\,\frac{dv_t}{dt}\bigg|_{t=0}\bigg\rangle_\omega\,\frac{\omega^n}{n!} + \sum\limits_{ |I|=p \atop |J|=q} u_{I\bar{J}}\,\overline{v_{I\bar{J}}(0)}\,\bigg(\sum\limits_{j=1}^n\gamma_j - \sum\limits_{l\in I}\gamma_l - \sum\limits_{r\in J}\gamma_r\bigg)\,\frac{\omega^n}{n!} \\
  & = & \bigg\langle u,\,\frac{dv_t}{dt}\bigg|_{t=0} + [\Lambda_\omega,\,\gamma\wedge\cdot]\,v_0\bigg\rangle_\omega\,\frac{\omega^n}{n!} = u\wedge\star_\omega\overline{\bigg(\frac{dv_t}{dt}\bigg|_{t=0} + [\Lambda_\omega,\,\gamma\wedge\cdot]\,v_0\bigg)},\end{eqnarray*} where the second equality follows from the standard pointwise formula (\ref{eqn:coordinates-standard-formula-anticom}) for $[\Lambda_\omega,\,\gamma\wedge\cdot]$ (see e.g. [Dem97, VI, $\S.5.2$, Proposition 5.8]) that we now reproduce.

Whenever $\omega$ is a Hermitian metric and $\gamma$ is a real $(1,\,1)$-form that are diagonalised simultaneously at a given point $x$ by a choice of local holomorphic coordinates $(z_1,\dots , z_n)$ centred at $x$ such that (\ref{eqn:omega-gamma_diagonalisation}) holds, for any bidegree $(p,\,q)$ and any $(p,\,q)$-form $v$ written in these local coordinates as $v=\sum_{|I|=p,\, |J|=q}v_{I\bar{J}}\,dz_I\wedge d\bar{z}_J$, the following equality holds at $x$: \begin{eqnarray}\label{eqn:coordinates-standard-formula-anticom} [\Lambda_\omega,\,\gamma\wedge\cdot]\,v = \sum\limits_{ |I|=p \atop |J|=q} v_{I\bar{J}}\,\bigg(\sum\limits_{j=1}^n\gamma_j - \sum\limits_{l\in I}\gamma_l - \sum\limits_{r\in J}\gamma_r\bigg)\,dz_I\wedge d\bar{z}_J.\end{eqnarray} In the special case where $\gamma=\omega$, we have $\gamma_1=\dots = \gamma_n=1$ and this standard formula reads: \begin{equation}\label{eqn:L-Lambda_comm_standard-formula}[\Lambda_\omega,\,L_\omega] = (n-p-q)\,\mbox{Id} \hspace{3ex}\mbox{on}\hspace{1ex} (p,\,q)-\mbox{forms}.\end{equation}

Since the form $u\in C^\infty_{p,\,q}(X,\,\C)$ was arbitrary, formula (\ref{eqn:1st-variation_Hodge-star}) follows.  \hfill $\Box$

\subsection{Computation of the first variation of the trace of $(p,\,q)$-forms}\label{subsection:1st-variation_trace_pq}

In this subsection, we generalise the result of $\S$\ref{subsection:1st-variation_trace_11} to the case of forms of arbitrary bidegree $(p,\,q)$ by using the first variation of the Hodge star operator computed in $\S$\ref{subsection:1st-variation_Hodge-star}.

We will need the following preliminary result.

\begin{Lem}\label{Lem:star-eta-wedge_commutation} Let $X$ be an $n$-dimensional complex manifold equipped with a Hermitian metric $\omega$ and let $\eta$ be a $(1,\,1)$-form on $X$. The following formula holds for operators acting on differential forms of any bidegree on $X$: \begin{equation}\label{eqn:star-eta-wedge_commutation}\star_\omega(\eta\wedge\cdot) = (\bar\eta\wedge\cdot)^\star_\omega\,\star_\omega,\end{equation} where $\eta\wedge\cdot$ is the operator of multiplication by $\eta$ and $(\bar\eta\wedge\cdot)^\star_\omega$ is its adjoint w.r.t. the pointwise inner product $\langle\,\cdot\,,\,\cdot\,\rangle_\omega$ defined by $\omega$.

\end{Lem}

\vspace{2ex}

An immediate consequence obtained when $\eta=\omega$ is the following well-known formula

\begin{Cor}\label{Cor:star-L-Lambda_commutation} Under the hypotheses of Lemma \ref{Lem:star-eta-wedge_commutation}, we have: \begin{equation}\label{eqn:star-L-Lambda_commutation}\star_\omega L_\omega = \Lambda_\omega\,\star_\omega,\end{equation} where $L_\omega:= \omega\wedge\cdot$ is the Lefschetz operator of multiplication by $\omega$ and $\Lambda_\omega = L_\omega^\star$ is its adjoint w.r.t. the pointwise inner product $\langle\,\cdot\,,\,\cdot\,\rangle_\omega$ defined by $\omega$.

\end{Cor}

\vspace{2ex}

\noindent {\it Proof of Lemma \ref{Lem:star-eta-wedge_commutation}.} Let $(p,\,q)$ be an arbitrary bidegree and let $u\in\Lambda^{p+1,\,q+1}T^\star X$, $v\in\Lambda^{p,\,q}T^\star X$ be arbitrary forms of the indicated bidegrees.

The definition of the Hodge star operator $\star_\omega$ yields the first and the third equalities below:  \begin{eqnarray}\label{eqn:star-L-Lambda_commutation_proof_1} u\wedge\star_\omega\overline{(\eta\wedge v)} & = & \langle u,\,\eta\wedge v\rangle_\omega\,dV_\omega = \langle (\eta\wedge\cdot)_\omega^\star u,\, v\rangle_\omega\,dV_\omega = (\eta\wedge\cdot)_\omega^\star u\wedge\star_\omega\bar{v}.\end{eqnarray} 

The formula to prove being pointiwse, we fix an arbitrary point $x\in X$ and we choose local holomorphic coordinates $(z_1,\dots , z_n)$ about $x$ such that $\omega$ is given by the identity matrix at $x$ in these coordinates. In particular, the adjoints of the multiplication operators $dz_j\wedge\cdot$ and $d\bar{z}_j\wedge\cdot$ w.r.t. the pointwise inner product induced by $\omega$ at $x$ are given by the contractions with the corresponding tangent vectors at $x$: \begin{eqnarray}\label{eqn:adjoints_dz_j_bar}(dz_j\wedge\cdot)^\star_\omega = \frac{\partial}{\partial z_j}\lrcorner\cdot, \hspace{5ex} (d\bar{z}_j\wedge\cdot)^\star_\omega = \frac{\partial}{\partial\bar{z}_j}\lrcorner\cdot,\end{eqnarray} for every $j$. Let $$\eta = \sum\limits_{j,\,k=1}^n\eta_{j\bar{k}}\,idz_j\wedge d\bar{z}_k$$ be the local expression of $\eta$.

The last term in (\ref{eqn:star-L-Lambda_commutation_proof_1}) reads: \begin{eqnarray}\label{eqn:star-L-Lambda_commutation_proof_2}(\eta\wedge\cdot)_\omega^\star u\wedge\star_\omega\bar{v} & = & - \sum\limits_{j,\,k=1}^n\bar\eta_{j\bar{k}}\,i\,\bigg(\frac{\partial}{\partial\bar{z}_k}\lrcorner\frac{\partial}{\partial z_j}\lrcorner u\bigg)\wedge\star_\omega\bar{v}.\end{eqnarray}

Meanwhile, we have:

\begin{eqnarray*}\bigg(\frac{\partial}{\partial\bar{z}_k}\lrcorner\frac{\partial}{\partial z_j}\lrcorner u\bigg)\wedge\star_\omega\bar{v} & = & \frac{\partial}{\partial\bar{z}_k}\lrcorner\bigg[\bigg(\frac{\partial}{\partial z_j}\lrcorner u\bigg)\wedge\star_\omega\bar{v}\bigg] - (-1)^{p+q+1}\,\bigg(\frac{\partial}{\partial z_j}\lrcorner u\bigg)\wedge\bigg(\frac{\partial}{\partial\bar{z}_k}\lrcorner\star_\omega\bar{v}\bigg) \\
  \nonumber & = & (-1)^{p+q}\,\bigg(\frac{\partial}{\partial z_j}\lrcorner u\bigg)\wedge\bigg(\frac{\partial}{\partial\bar{z}_k}\lrcorner\star_\omega\bar{v}\bigg) \\
  \nonumber & = & (-1)^{p+q}\,\frac{\partial}{\partial z_j}\lrcorner\bigg[u\wedge\bigg(\frac{\partial}{\partial\bar{z}_k}\lrcorner\star_\omega\bar{v}\bigg)\bigg] - u\wedge\bigg(\frac{\partial}{\partial z_j}\lrcorner\frac{\partial}{\partial\bar{z}_k}\lrcorner\star_\omega\bar{v}\bigg) \\
  & = & - u\wedge\bigg(\frac{\partial}{\partial z_j}\lrcorner\frac{\partial}{\partial\bar{z}_k}\lrcorner\star_\omega\bar{v}\bigg),\end{eqnarray*} where the first term on the r.h.s. of each of the first and third lines above vanishes for bidegree reasons. Indeed, $((\partial/\partial z_j)\lrcorner u)\wedge\star_\omega\bar{v}$ is of bidegree $(n,\,n+1)$ and $u\wedge((\partial/\partial\bar{z}_k)\lrcorner\star_\omega\bar{v})$ is of bidegree $(n+1,\,n)$ since $u$ is of bidegree $(p+1,\,q+1)$ and $\star_\omega\bar{v}$ is of bidegree $(n-p,\,n-q)$. Using again (\ref{eqn:adjoints_dz_j_bar}), this translates to

\begin{eqnarray}\label{eqn:star-L-Lambda_commutation_proof_3}\bigg(\frac{\partial}{\partial\bar{z}_k}\lrcorner\frac{\partial}{\partial z_j}\lrcorner u\bigg)\wedge\star_\omega\bar{v} = - u\wedge(d\bar{z}_k\wedge dz_j\wedge\cdot)^\star_\omega\,(\star_\omega\bar{v}).\end{eqnarray}

Putting (\ref{eqn:star-L-Lambda_commutation_proof_2}) and (\ref{eqn:star-L-Lambda_commutation_proof_3}) together, we get:

\begin{eqnarray}\label{eqn:star-L-Lambda_commutation_proof_4}\nonumber(\eta\wedge\cdot)_\omega^\star u\wedge\star_\omega\bar{v} & = & \sum\limits_{j,k=1}^n\bar\eta_{j\bar{k}}\,i\, u\wedge(d\bar{z}_k\wedge dz_j\wedge\cdot)^\star_\omega\,(\star_\omega\bar{v}) = -u\wedge\bigg(\sum\limits_{j,\,k=1}^n\eta_{j\bar{k}}\,i\,d\bar{z}_k\wedge dz_j\wedge\cdot\bigg)^\star_\omega\,(\star_\omega\bar{v}) \\
 & = & u\wedge(\eta\wedge\cdot)^\star_\omega\,(\star_\omega\bar{v}).\end{eqnarray}

Finally, putting (\ref{eqn:star-L-Lambda_commutation_proof_1}) and (\ref{eqn:star-L-Lambda_commutation_proof_4}) together, we get: \begin{eqnarray*}u\wedge\star_\omega\overline{(\eta\wedge v)} = u\wedge\overline{(\bar\eta\wedge\cdot)^\star_\omega\,(\star_\omega v)}\end{eqnarray*} for all forms $u\in\Lambda^{p+1,\,q+1}T^\star X$, $v\in\Lambda^{p,\,q}T^\star X$. This proves formula (\ref{eqn:star-eta-wedge_commutation}).  \hfill $\Box$

\vspace{2ex}

We are now in a position to prove the following generalisation of Lemma \ref{Lem:1st-variation_trace_11} to the case of forms of an arbitrary bidegree $(p,\,q)$.

\begin{Lem}\label{Lem:1st-variation_trace_pq} For any bidegree $(p,\,q)$ and any $C^\infty$ family $(\alpha_t)_{t\in(-\varepsilon,\,\varepsilon)}$ of forms $\alpha_t\in C^\infty_{p,\,q}(X,\,\C)$ with $\varepsilon>0$ so small that $\omega+t\gamma>0$ for all $t\in(-\varepsilon,\,\varepsilon)$, the following formula holds: \begin{eqnarray}\label{eqn:1st-variation_trace_pq}\frac{d}{dt}\bigg|_{t=0}\,(\Lambda_{\omega+t\gamma}\alpha_t) = \Lambda_\omega\bigg(\frac{d\alpha_t}{dt}\bigg|_{t=0}\bigg) - (\gamma\wedge\cdot)^\star_\omega\,\alpha_0,  \hspace{10ex} t\in(-\varepsilon,\,\varepsilon).\end{eqnarray}
     
\end{Lem}

\noindent {\it Proof.} From formula (\ref{eqn:star-L-Lambda_commutation}) applied to the metric $\omega+t\gamma$ and from $\star_{\omega+t\gamma}\,\star_{\omega+t\gamma} = (-1)^{p+q}\,\mbox{Id}$ on $(p,\,q)$-forms, we get: \begin{eqnarray*}\Lambda_{\omega+t\gamma}\alpha_t = (-1)^{p+q}\,\star_{\omega+t\gamma}\,L_{\omega+t\gamma}\,\star_{\omega+t\gamma}\alpha_t = (-1)^{p+q}\,\star_{\omega+t\gamma}\,v_t,\end{eqnarray*} where we put $v_t:=L_{\omega+t\gamma}(\star_{\omega+t\gamma}\alpha_t)$.

Applying formula (\ref{eqn:1st-variation_Hodge-star}), we get the second equality below: \begin{eqnarray}\label{eqn:1st-variation_trace_pq_proof_1}\frac{d}{dt}\bigg|_{t=0}\,(\Lambda_{\omega+t\gamma}\alpha_t) = (-1)^{p+q}\,\frac{d}{dt}\bigg|_{t=0}\,(\star_{\omega+t\gamma}\,v_t) =  (-1)^{p+q}\,\star_\omega\bigg(\frac{dv_t}{dt}\bigg|_{t=0} + [\Lambda_\omega,\,\gamma\wedge\cdot]\,v_0\bigg).\end{eqnarray}

On the other hand, the same formula (\ref{eqn:1st-variation_Hodge-star}) yields the third equality below: \begin{eqnarray}\label{eqn:1st-variation_trace_pq_proof_2}\nonumber\frac{dv_t}{dt}\bigg|_{t=0} & = & \frac{d}{dt}\bigg|_{t=0}\bigg((\omega+t\gamma)\wedge\star_{\omega+t\gamma}\alpha_t\bigg) = \gamma\wedge\star_\omega\alpha_0 + \omega\wedge \frac{d}{dt}\bigg|_{t=0}\bigg(\star_{\omega+t\gamma}\alpha_t\bigg) \\
  & = & \gamma\wedge\star_\omega\alpha_0 + \omega\wedge\star_\omega\bigg(\frac{d\alpha_t}{dt}\bigg|_{t=0} + [\Lambda_\omega,\,\gamma\wedge\cdot]\,\alpha_0\bigg).\end{eqnarray}

Putting (\ref{eqn:1st-variation_trace_pq_proof_1}) and (\ref{eqn:1st-variation_trace_pq_proof_2}) together, we get: \begin{eqnarray}\label{eqn:1st-variation_trace_pq_proof_3}\nonumber\frac{d}{dt}\bigg|_{t=0}\,(\Lambda_{\omega+t\gamma}\alpha_t) & = & (-1)^{p+q}\,\star_\omega(\gamma\wedge\star_\omega\alpha_0) + \Lambda_\omega\bigg(\frac{d\alpha_t}{dt}\bigg|_{t=0} + [\Lambda_\omega,\,\gamma\wedge\cdot]\,\alpha_0\bigg) \\
  & + & (-1)^{p+q}\,\star_\omega\bigg([\Lambda_\omega,\,\gamma\wedge\cdot]\,L_\omega(\star_\omega\alpha_0)\bigg),\end{eqnarray} where the general equality $\star_\omega(\omega\wedge\cdot) = \Lambda_\omega\,\star_\omega$ was used to get the second term on the r.h.s. and the equality $v_0 = L_\omega(\star_\omega\alpha_0)$ lead to the third term.

To compute the last term in (\ref{eqn:1st-variation_trace_pq_proof_3}), we notice that  \begin{eqnarray}\label{eqn:1st-variation_trace_pq_proof_4}\nonumber[\Lambda_\omega,\,\gamma\wedge\cdot]\,L_\omega(\star_\omega\alpha_0) & = & \Lambda_\omega(\omega\wedge\star_\omega\alpha_0\wedge\gamma) - \Lambda_\omega(\omega\wedge\star_\omega\alpha_0)\wedge\gamma \\
  \nonumber & = & [\Lambda_\omega,\,L_\omega]\,(\gamma\wedge\star_\omega\alpha_0) + \omega\wedge\Lambda_\omega(\gamma\wedge\star_\omega\alpha_0) - [\Lambda_\omega,\,L_\omega]\,(\star_\omega\alpha_0)\wedge\gamma - \omega\wedge\gamma\wedge\Lambda_\omega(\star_\omega\alpha_0) \\
  \nonumber & = & (p+q-n-2)\,\gamma\wedge\star_\omega\alpha_0  - (p+q-n)\,\gamma\wedge\star_\omega\alpha_0 + \omega\wedge[\Lambda_\omega,\,\gamma\wedge\cdot]\,(\star_\omega\alpha_0) \\
  \nonumber & = & -2\,\gamma\wedge\star_\omega\alpha_0 + \omega\wedge[\Lambda_\omega,\,\gamma\wedge\cdot]\,(\star_\omega\alpha_0),\end{eqnarray} where we have used the standard formula (\ref{eqn:L-Lambda_comm_standard-formula}). Thus, (\ref{eqn:1st-variation_trace_pq_proof_3}) becomes: \begin{eqnarray}\label{eqn:1st-variation_trace_pq_proof_5}\nonumber\frac{d}{dt}\bigg|_{t=0}\,(\Lambda_{\omega+t\gamma}\alpha_t) & = & \Lambda_\omega\bigg(\frac{d\alpha_t}{dt}\bigg|_{t=0}\bigg) + \Lambda_\omega\bigg([\Lambda_\omega,\,\gamma\wedge\cdot]\,\alpha_0\bigg) - (\gamma\wedge\cdot)^\star_\omega\,\alpha_0 \\
  & + & (-1)^{p+q}\,\Lambda_\omega\bigg(\star_\omega\,[\Lambda_\omega,\,\gamma\wedge\cdot]\,\star_\omega\,\alpha_0\bigg),\end{eqnarray} where we have used formula (\ref{eqn:star-eta-wedge_commutation}) (with $\eta=\gamma$, a {\it real} form) and the fact that $\star_\omega\star_\omega\alpha_0 = (-1)^{p+q}\,\alpha_0$ to get $\star_\omega\,(\gamma\wedge\star_\omega\alpha_0) = (-1)^{p+q}\,(\gamma\wedge\cdot)^\star_\omega\,\alpha_0$.

Next, to compute the last term in (\ref{eqn:1st-variation_trace_pq_proof_5}), we notice the following identities:  \begin{eqnarray}\label{eqn:1st-variation_trace_pq_proof_6}\nonumber \star_\omega\,[\Lambda_\omega,\,\gamma\wedge\cdot]\,\star_\omega\,\alpha_0 & = & \star_\omega\,\Lambda_\omega\, (\gamma\wedge\cdot)\,\star_\omega\,\alpha_0 - \star_\omega\,(\gamma\wedge\cdot)\,\Lambda_\omega\,\star_\omega\,\alpha_0 \\
  \nonumber & = & L_\omega\,\bigg(\star_\omega\,(\gamma\wedge\cdot)\bigg)\,\star_\omega\,\alpha_0 - (\gamma\wedge\cdot)^\star_\omega\,\bigg(\star_\omega\,\Lambda_\omega\bigg)\,\star_\omega\,\alpha_0  \\
  \nonumber & = & (-1)^{p+q}\, L_\omega\,(\gamma\wedge\cdot)^\star_\omega\,\alpha_0 -  (-1)^{p+q}\,(\gamma\wedge\cdot)^\star_\omega\,L_\omega\,\alpha_0 \\
        & = & (-1)^{p+q}\,[L_\omega,\,(\gamma\wedge\cdot)^\star_\omega]\,\alpha_0 = -(-1)^{p+q}\,[\Lambda_\omega,\,\gamma\wedge\cdot]^\star_\omega\,\alpha_0,\end{eqnarray} where we have used several times the identities $\star_\omega\,(\gamma\wedge\cdot) = (\gamma\wedge\cdot)^\star_\omega\,\star_\omega$ (see formula (\ref{eqn:star-eta-wedge_commutation})) and $\star_\omega\,\Lambda_\omega = L_\omega\,\star_\omega$ (see formula (\ref{eqn:star-L-Lambda_commutation})).

Putting (\ref{eqn:1st-variation_trace_pq_proof_5}) and (\ref{eqn:1st-variation_trace_pq_proof_6}) together, we get:

\begin{eqnarray}\label{eqn:1st-variation_trace_pq_proof_7}\frac{d}{dt}\bigg|_{t=0}\,(\Lambda_{\omega+t\gamma}\alpha_t) & = & \Lambda_\omega\bigg(\frac{d\alpha_t}{dt}\bigg|_{t=0}\bigg) - (\gamma\wedge\cdot)^\star_\omega\,\alpha_0 + \Lambda_\omega\bigg([\Lambda_\omega,\,\gamma\wedge\cdot] - [\Lambda_\omega,\,\gamma\wedge\cdot]^\star_\omega\bigg)\,\alpha_0.\end{eqnarray} Now, the operator $[\Lambda_\omega,\,\gamma\wedge\cdot]$ is of order zero, so it acts pointwise on differential forms. The standard formula (\ref{eqn:coordinates-standard-formula-anticom}) and the fact that the eigenvalues $\gamma_1,\dots , \gamma_n$ of $\gamma$ w.r.t. $\omega$ are {\it real} (since $\gamma$ is) imply at once that the operator $[\Lambda_\omega,\,\gamma\wedge\cdot]$ is self-adjoint w.r.t. the pointwise inner product defined by $\omega$. Therefore, the last term in (\ref{eqn:1st-variation_trace_pq_proof_7}) vanishes, so (\ref{eqn:1st-variation_trace_pq_proof_7}) amounts to (\ref{eqn:1st-variation_trace_pq}). This completes the proof of Lemma \ref{Lem:1st-variation_trace_pq}. \hfill $\Box$

\vspace{2ex}

Note that, when $(p,\,q)=(1,\,1)$, an immediate computation in coordinates yields: \begin{eqnarray}\label{eqn:mult-adjoint-inner-prod_11}(\eta\wedge\cdot)^\star_\omega = \langle\cdot,\,\eta\rangle_\omega  \hspace{3ex} \mbox{on}\hspace{1ex} \Lambda^{1,\,1}T^\star X,\end{eqnarray} so Lemma \ref{Lem:1st-variation_trace_pq} reproves Lemma \ref{Lem:1st-variation_trace_11}.

\subsection{Computation of the first variation of the adjoints and the Laplacians}\label{subsection:1st-variation_adjoints-Laplacians}

Consider the setting of Lemma \ref{Lem:1st-variation_Hodge-star}. We start by deriving the standard formula $\partial^\star_{\omega+t\gamma}v_t = -\star_{\omega+t\gamma}\bar\partial\star_{\omega+t\gamma}v_t = -\star_{\omega+t\gamma}\alpha_t$, where we put $\alpha_t:=\bar\partial\star_{\omega+t\gamma}v_t$. Applying formula (\ref{eqn:1st-variation_Hodge-star}) with $\alpha_t$ in place of $v_t$, we get the second equality below: \begin{eqnarray*}\frac{d}{dt}\bigg|_{t=0}(\partial^\star_{\omega+t\gamma}v_t) & = & -\frac{d}{dt}\bigg|_{t=0}(\star_{\omega+t\gamma}\alpha_t) = -\star_\omega\bigg(\frac{d\alpha_t}{dt}\bigg|_{t=0} + [\Lambda_\omega,\,\gamma\wedge\cdot]\,\alpha_0\bigg) \\
  & = & - \star_\omega\bar\partial\bigg(\frac{d}{dt}\bigg|_{t=0}(\star_{\omega+t\gamma}v_t)\bigg) - \star_\omega\bigg([\Lambda_\omega,\,\gamma\wedge\cdot]\,\bar\partial\star_\omega v_0\bigg).\end{eqnarray*} Applying again formula (\ref{eqn:1st-variation_Hodge-star}) to the first term on the last line above, using the standard formula $\partial^\star_\omega = -\star_\omega\bar\partial\star_\omega$ and the fact that $\star_\omega\star_\omega = (-1)^k\,\mbox{Id}$ on $k$-forms (for any $k$), we get formula (\ref{eqn:1st-variation_del-adjoint}) in

\begin{Lem}\label{Lem:1st-variation_adjoints} (i)\, For every bidegree $(p,\,q)$ and every $C^\infty$ family $(v_t)_{t\in(-\varepsilon,\,\varepsilon)}$ of forms $v_t\in C^\infty_{p,\,q}(X,\,\C)$ with $\varepsilon>0$ so small that $\omega+t\gamma>0$ for all $t\in(-\varepsilon,\,\varepsilon)$, the following formulae hold: \begin{eqnarray}\label{eqn:1st-variation_del-adjoint}\frac{d}{dt}\bigg|_{t=0}(\partial^\star_{\omega+t\gamma}v_t) = \partial^\star_\omega\bigg(\frac{dv_t}{dt}\bigg|_{t=0}\bigg) + \partial^\star_\omega\bigg([\Lambda_\omega,\,\gamma\wedge\cdot]\,v_0\bigg) + (-1)^{\deg v_0 +1}
\bigg(\star_\omega[\Lambda_\omega,\,\gamma\wedge\cdot]\star_\omega\bigg)\partial^\star_\omega v_0,\end{eqnarray}  \begin{eqnarray}\label{eqn:1st-variation_del-bar-adjoint}\frac{d}{dt}\bigg|_{t=0}(\bar\partial^\star_{\omega+t\gamma}v_t) = \bar\partial^\star_\omega\bigg(\frac{dv_t}{dt}\bigg|_{t=0}\bigg) + \bar\partial^\star_\omega\bigg([\Lambda_\omega,\,\gamma\wedge\cdot]\,v_0\bigg) + (-1)^{\deg v_0 +1}
    \bigg(\star_\omega[\Lambda_\omega,\,\gamma\wedge\cdot]\star_\omega\bigg)\bar\partial^\star_\omega v_0.\end{eqnarray}

  \vspace{1ex}

 (ii)\, For every degree $k$ and every $C^\infty$ family $(v_t)_{t\in(-\varepsilon,\,\varepsilon)}$ of forms $v_t\in C^\infty_k(X,\,\C)$ with $\varepsilon>0$ so small that $\omega+t\gamma>0$ for all $t\in(-\varepsilon,\,\varepsilon)$, the following formula holds:   \begin{eqnarray}\label{eqn:1st-variation_d-adjoint}\frac{d}{dt}\bigg|_{t=0}(d^\star_{\omega+t\gamma}v_t) = d^\star_\omega\bigg(\frac{dv_t}{dt}\bigg|_{t=0}\bigg) + d^\star_\omega\bigg([\Lambda_\omega,\,\gamma\wedge\cdot]\,v_0\bigg) + (-1)^{\deg v_0 +1}
\bigg(\star_\omega[\Lambda_\omega,\,\gamma\wedge\cdot]\star_\omega\bigg)d^\star_\omega v_0.\end{eqnarray}

\end{Lem}

\noindent {\it Proof.} Formula (\ref{eqn:1st-variation_del-adjoint}) was proved just before the statement, (\ref{eqn:1st-variation_del-bar-adjoint}) follows by conjugating (\ref{eqn:1st-variation_del-adjoint}), while (\ref{eqn:1st-variation_d-adjoint}) follows by adding (\ref{eqn:1st-variation_del-adjoint}) and (\ref{eqn:1st-variation_del-bar-adjoint}) up for each of the pure-type components of $v_t$.  \hfill $\Box$

\vspace{3ex}

Let us now consider the $\partial$-Laplacians induced by $\omega$ and $\omega+t\gamma$ (with $t\in\R$ close to $0$): \begin{eqnarray*}\Delta'_\omega = \partial\partial^\star_\omega + \partial^\star_\omega\partial, \hspace{6ex} \Delta'_{\omega+t\gamma} = \partial\partial^\star_{\omega+t\gamma} + \partial^\star_{\omega+t\gamma}\partial\end{eqnarray*} and their $\bar\partial$- and $d$-analogues $\Delta''_\omega$, $\Delta''_{\omega+t\gamma}$, $\Delta_\omega$, $\Delta_{\omega+t\gamma}$. We get:

\begin{eqnarray*}\frac{d}{dt}\bigg|_{t=0}(\Delta'_{\omega+t\gamma}v_t) = \partial\bigg(\frac{d}{dt}\bigg|_{t=0}(\partial^\star_{\omega+t\gamma}v_t)\bigg) + \frac{d}{dt}\bigg|_{t=0}(\partial^\star_{\omega+t\gamma}(\partial v_t))\end{eqnarray*} and its analogues for $\Delta''_\omega$ and $\Delta_\omega$, $\Delta_{\omega+t\gamma}$.

Then, a straightforward application of formulae (\ref{eqn:1st-variation_del-adjoint})--(\ref{eqn:1st-variation_d-adjoint}) of Lemma \ref{Lem:1st-variation_adjoints} yields the following

  \begin{Lem}\label{Lem:1st-variation_Laplacians} (i)\, For every bidegree $(p,\,q)$ and every $C^\infty$ family $(v_t)_{t\in(-\varepsilon,\,\varepsilon)}$ of forms $v_t\in C^\infty_{p,\,q}(X,\,\C)$ with $\varepsilon>0$ so small that $\omega+t\gamma>0$ for all $t\in(-\varepsilon,\,\varepsilon)$, the following formulae hold:

    \begin{eqnarray}\label{eqn:1st-variation_del-Laplacian}\frac{d}{dt}\bigg|_{t=0}(\Delta'_{\omega+t\gamma}v_t) & = & \Delta'_\omega\bigg(\frac{dv_t}{dt}\bigg|_{t=0}\bigg) + \partial\partial^\star_\omega\bigg([\Lambda_\omega,\,\gamma\wedge\cdot]\,v_0\bigg) + \partial^\star_\omega\bigg([\Lambda_\omega,\,\gamma\wedge\cdot]\,\partial v_0\bigg) \\
      \nonumber      & + & (-1)^{\deg v_0 +1}\,\partial\bigg(\star_\omega[\Lambda_\omega,\,\gamma\wedge\cdot]\star_\omega\bigg)\partial^\star_\omega v_0 + (-1)^{\deg v_0}\,\bigg(\star_\omega[\Lambda_\omega,\,\gamma\wedge\cdot]\star_\omega\bigg)\,\partial^\star_\omega\partial v_0,\end{eqnarray}

    \begin{eqnarray}\label{eqn:1st-variation_del-bar-Laplacian}\frac{d}{dt}\bigg|_{t=0}(\Delta''_{\omega+t\gamma}v_t) & = & \Delta''_\omega\bigg(\frac{dv_t}{dt}\bigg|_{t=0}\bigg) + \bar\partial\bar\partial^\star_\omega\bigg([\Lambda_\omega,\,\gamma\wedge\cdot]\,v_0\bigg) + \bar\partial^\star_\omega\bigg([\Lambda_\omega,\,\gamma\wedge\cdot]\,\bar\partial v_0\bigg) \\
      \nonumber      & + & (-1)^{\deg v_0 +1}\,\bar\partial\bigg(\star_\omega[\Lambda_\omega,\,\gamma\wedge\cdot]\star_\omega\bigg)\bar\partial^\star_\omega v_0 + (-1)^{\deg v_0}\,\bigg(\star_\omega[\Lambda_\omega,\,\gamma\wedge\cdot]\star_\omega\bigg)\,\bar\partial^\star_\omega\bar\partial v_0.\end{eqnarray}

     \vspace{1ex}

   (ii)\, For every degree $k$ and every $C^\infty$ family $(v_t)_{t\in(-\varepsilon,\,\varepsilon)}$ of forms $v_t\in C^\infty_k(X,\,\C)$ with $\varepsilon>0$ so small that $\omega+t\gamma>0$ for all $t\in(-\varepsilon,\,\varepsilon)$, the following formula holds:

   \begin{eqnarray}\label{eqn:1st-variation_d-Laplacian}\frac{d}{dt}\bigg|_{t=0}(\Delta_{\omega+t\gamma}v_t) & = & \Delta_\omega\bigg(\frac{dv_t}{dt}\bigg|_{t=0}\bigg) + dd^\star_\omega\bigg([\Lambda_\omega,\,\gamma\wedge\cdot]\,v_0\bigg) + d^\star_\omega\bigg([\Lambda_\omega,\,\gamma\wedge\cdot]\,d v_0\bigg) \\
      \nonumber      & + & (-1)^{\deg v_0 +1}\,d\bigg(\star_\omega[\Lambda_\omega,\,\gamma\wedge\cdot]\star_\omega\bigg)d^\star_\omega v_0 + (-1)^{\deg v_0}\,\bigg(\star_\omega[\Lambda_\omega,\,\gamma\wedge\cdot]\star_\omega\bigg)\,d^\star_\omega d v_0.\end{eqnarray}

\end{Lem}

\subsection{Computation of the first variation of the orthogonal projections onto the kernels of the Laplacians}\label{subsection:1st-variation_projections-ker-Laplacians}

We continue with the setting of Lemma \ref{Lem:1st-variation_Hodge-star}. We will need a formula for the first variation of the $L^2_\omega$-orthogonal projection $P^{(3)}_\omega:C^\infty_3(X,\,\C)\longrightarrow\mbox{Im}\,d$ induced by the three-space decomposition (\ref{eqn:3-space-decomp_d-3}). Since $\partial\omega\in\ker d = {\cal H}^3_\Delta(X,\,\C)\oplus\mbox{Im}\,d$ under our SKT assumption, we will only be dealing with the restriction of $P^{(3)}_\omega$ to $\ker d$, so the computation of the first variation of $P^{(3)}_\omega$ is equivalent to the computation of the first variation of the $L^2_\omega$-orthogonal projection $F^{(3)}_\Delta:C^\infty_3(X,\,\C)\longrightarrow{\cal H}^3_\Delta(X,\,\C)$.

Allowing for more flexibility, we consider an arbitrary bidegree $(p,\,q)$ and the $L^2_{\omega+t\gamma}$-orthogonal projection \begin{equation*}F'_t:C^\infty_{p,\,q}(X,\,\C)\longrightarrow{\cal H}^{p,\,q}_{\Delta'_{\omega+t\gamma}}(X,\,\C), \hspace{5ex} t\sim 0,\end{equation*} onto the kernel ${\cal H}^{p,\,q}_{\Delta'_{\omega+t\gamma}}(X,\,\C)$ of the $\partial$-Laplacian $\Delta'_{\omega+t\gamma}$. This projection is induced by the three-space $L^2_{\omega+t\gamma}$-orthogonal decomposition  \begin{equation}\label{eqn:3-space-decomp_del-pq}C^\infty_{p,\,q}(X,\,\C) = {\cal H}^{p,\,q}_{\Delta'_{\omega+t\gamma}}(X,\,\C)\oplus\mbox{Im}\,\partial\oplus\mbox{Im}\,\partial^\star_{\omega+t\gamma}\end{equation} in which $\ker\partial = {\cal H}^{p,\,q}_{\Delta'_{\omega+t\gamma}}(X,\,\C)\oplus\mbox{Im}\,\partial$ and $\ker\partial^\star_{\omega+t\gamma} = {\cal H}^{p,\,q}_{\Delta'_{\omega+t\gamma}}(X,\,\C)\oplus\mbox{Im}\,\partial^\star_{\omega+t\gamma}$.

We will denote by $F''_t:C^\infty_{p,\,q}(X,\,\C)\longrightarrow{\cal H}^{p,\,q}_{\Delta''_{\omega+t\gamma}}(X,\,\C)$ and $F_t:C^\infty_k(X,\,\C)\longrightarrow{\cal H}^k_{\Delta_{\omega+t\gamma}}(X,\,\C)$ the analogous orthogonal projections onto the kernels of $\Delta''_{\omega+t\gamma}$ and $\Delta_{\omega+t\gamma}$.

Note that $(\omega+t\gamma)_{|t|<\varepsilon}$ is a $C^\infty$ family of Hermitian metrics on $X$, so the family $(\Delta'_{\omega+t\gamma})_{|t|<\varepsilon}$ is a $C^\infty$ family of elliptic operators on $C^\infty_{p,\,q}(X,\,\C)$. Moreover, we have Hodge isomorphisms: \begin{equation*}\ker\Delta'_{\omega+t\gamma}\simeq H^{p,\,q}_\partial(X,\,\C),  \hspace{6ex} t\in(-\varepsilon,\,\varepsilon),\end{equation*} so the dimension of the kernel of $\Delta'_{\omega+t\gamma}$ is independent of $t$. By the classical Kodaira-Spencer Theorem 5 in [KS60], the operator $F'_t$ depends in a $C^\infty$ way on $t\in(-\varepsilon,\,\varepsilon)$. 

We now set out to compute the first variation at $t=0$ of $F'_t$ using the classical Kodaira-Spencer {\it Cauchy-type formula} (see [KS60, p. 56] or [Kod86, chapter 7, (7.42)]):  \begin{equation}\label{eqn:K-S_Cauchy}F'_tv = -\frac{1}{2\pi i}\,\int_{\zeta\in C}(\Delta'_{\omega+t\gamma}-\zeta)^{-1}v\,d\zeta, \hspace{5ex} v\in C^\infty_{p,\,q}(X,\,\C),\, t\in(-\varepsilon,\,\varepsilon),\end{equation} where $C\subset\C$ is any Jordan curve in the complex plane that does not meet the spectrum of $\Delta'_\omega$. Recall that $\mbox{Spec}\,(\Delta'_\omega)$ is discrete and has $+\infty$ as its only accumulation point thanks to $\Delta'_\omega$ being elliptic and to $X$ being compact. Moreover, when $\varepsilon>0$ is small enough, $C\subset\C$ does not meet the spectrum of $\Delta'_{\omega+t\gamma}$ for any  $t\in(-\varepsilon,\,\varepsilon)$. (See [Kod86, chapter 7, (7.41)].)

In our case, we choose $C$ to be a circle centred at the origin in $\C$ so small that $0$ is the only eigenvalue of $\Delta'_\omega$ lying in the interior of $C$
 and all the other eigenvalues of $\Delta'_\omega$ lie outside the closed disc $\overline{D}$ bounded by $C$. Since $\mbox{dim}\ker\Delta'_{\omega+t\gamma}$ is independent of $t$, $0$ continues to be the only eigenvalue of $\Delta'_{\omega+t\gamma}$ lying in the interior of $C$ with all the other eigenvalues lying outside $\overline{D}$ for all $t\in(-\varepsilon,\,\varepsilon)$ if $\varepsilon>0$ is small enough.

Deriving (\ref{eqn:K-S_Cauchy}) w.r.t. $t$ and then taking $t=0$, we get \begin{equation}\label{eqn:K-S_Cauchy_derived_1}\frac{d}{dt}\bigg|_{t=0}(F'_tv) = -\frac{1}{2\pi i}\,\int_{\zeta\in C}\frac{d}{dt}\bigg|_{t=0}\bigg((\Delta'_{\omega+t\gamma}-\zeta)^{-1}v\bigg)\,d\zeta, \hspace{5ex} v\in C^\infty_{p,\,q}(X,\,\C),\, t\in(-\varepsilon,\,\varepsilon).\end{equation}

To continue the computations, we fix a form $v\in C^\infty_{p,\,q}(X,\,\C)$ and we put $$u_{t,\,\zeta}:=(\Delta'_{\omega+t\gamma}-\zeta)^{-1}v,  \hspace{6ex} t\in(-\varepsilon,\,\varepsilon), \, \zeta\in D_1,$$ where $D_1$ is an open disc centred at the origin in $\C$ that contains $\overline{D}$ such that all the positive eigenvalues of $\Delta'_{\omega+t\gamma}$ lie outside $D_1$ for all $t\in(-\varepsilon,\,\varepsilon)$. Thus, $$v= (\Delta'_{\omega+t\gamma}-\zeta)\, u_{t,\,\zeta},  \hspace{6ex} t\in(-\varepsilon,\,\varepsilon), \, \zeta\in D_1.$$

Since $v$ is independent of $t$, we get the first equality below: \begin{eqnarray*}0 = \frac{dv}{dt}\bigg|_{t=0} = \frac{d}{dt}\bigg|_{t=0}\bigg(\Delta'_{\omega+t\gamma}u_{t,\,\zeta}\bigg) - \zeta\, \frac{du_{t,\,\zeta}}{dt}\bigg|_{t=0} =  (\Delta'_\omega-\zeta)\,\bigg(\frac{du_{t,\,\zeta}}{dt}\bigg|_{t=0}\bigg) - A(u_{0,\,\zeta}),\end{eqnarray*} where the last equality follows from formula (\ref{eqn:1st-variation_del-Laplacian}) after denoting, for every $\zeta\in D_1$: \begin{eqnarray}\label{eqn:A'u_zeta_notation}A'(u_{0,\,\zeta}) & := & -\partial\partial^\star_\omega\bigg([\Lambda_\omega,\,\gamma\wedge\cdot]\,u_{0,\,\zeta}\bigg) - \partial^\star_\omega\bigg([\Lambda_\omega,\,\gamma\wedge\cdot]\,\partial u_{0,\,\zeta}\bigg) \\
  \nonumber      & + & (-1)^{\deg u_{0,\,\zeta}}\,\partial\bigg(\star_\omega[\Lambda_\omega,\,\gamma\wedge\cdot]\star_\omega\bigg)\partial^\star_\omega u_{0,\,\zeta} + (-1)^{\deg u_{0,\,\zeta}+1}\,\bigg(\star_\omega[\Lambda_\omega,\,\gamma\wedge\cdot]\star_\omega\bigg)\,\partial^\star_\omega\partial u_{0,\,\zeta}.\end{eqnarray}

We infer that, for every $v\in C^\infty_{p,\,q}(X,\,\C)$ and every $t\in(-\varepsilon,\,\varepsilon)$, we have: \begin{equation}\label{eqn:K-S_Cauchy_derived_2}\frac{d}{dt}\bigg|_{t=0}\bigg((\Delta'_{\omega+t\gamma}-\zeta)^{-1}v\bigg) = \frac{du_{t,\,\zeta}}{dt}\bigg|_{t=0} = (\Delta'_\omega-\zeta)^{-1}\,\bigg(A'(u_{0,\,\zeta})\bigg).\end{equation}

From this, we easily get the following

\begin{Lem}\label{Lem:1st-variation_projections-ker-Laplacians} (i)\, For every bidegree $(p,\,q)$ and every form $v\in C^\infty_{p,\,q}(X,\,\C)$, the following formulae hold: \begin{eqnarray}\label{eqn:1st-variation_projection-ker-Delta'-Delta''}\frac{d}{dt}\bigg|_{t=0}(F'_tv) = F'_0(A'(u_{0,\,0})) \hspace{3ex} \mbox{and} \hspace{3ex} \frac{d}{dt}\bigg|_{t=0}(F''_tv) = F''_0(A''(u_{0,\,0})),\end{eqnarray} where $A'(u_{0,\,0})$ is the form defined in (\ref{eqn:A'u_zeta_notation}) for $\zeta=0$ and $A''(u_{0,\,0})$ is the analogous form defined by replacing $\partial$ with $\bar\partial$.

  \vspace{1ex}

   (ii)\, For every degree $k$ and every $C^\infty$ family $(v_t)_{t\in(-\varepsilon,\,\varepsilon)}$ of forms $v_t\in C^\infty_k(X,\,\C)$ with $\varepsilon>0$ so small that $\omega+t\gamma>0$ for all $t\in(-\varepsilon,\,\varepsilon)$, the following formula holds:  \begin{eqnarray}\label{eqn:1st-variation_projection-ker-Delta}\frac{d}{dt}\bigg|_{t=0}(F_tv) = F_0(A(u_{0,\,0})),\end{eqnarray} where, for every $\zeta\in D_1$, we put: \begin{eqnarray}\label{eqn:Au_zeta_notation}A(u_{0,\,\zeta}) & := & -dd^\star_\omega\bigg([\Lambda_\omega,\,\gamma\wedge\cdot]\,u_{0,\,\zeta}\bigg) - d^\star_\omega\bigg([\Lambda_\omega,\,\gamma\wedge\cdot]\,d u_{0,\,\zeta}\bigg) \\
  \nonumber      & + & (-1)^{\deg u_{0,\,\zeta}}\,d\bigg(\star_\omega[\Lambda_\omega,\,\gamma\wedge\cdot]\star_\omega\bigg)d^\star_\omega u_{0,\,\zeta} + (-1)^{\deg u_{0,\,\zeta}+1}\,\bigg(\star_\omega[\Lambda_\omega,\,\gamma\wedge\cdot]\star_\omega\bigg)\,d^\star_\omega d u_{0,\,\zeta}.\end{eqnarray}

\end{Lem}

\noindent {\it Proof.} It suffices to prove the first equality in (\ref{eqn:1st-variation_projection-ker-Delta'-Delta''}) since the second one is obtained from the first by taking conjugates and (\ref{eqn:1st-variation_projection-ker-Delta}) follows by running the same argument with $d$ in place of $\partial$ and using (\ref{eqn:1st-variation_d-Laplacian}) rather than (\ref{eqn:1st-variation_del-Laplacian}).

Putting (\ref{eqn:K-S_Cauchy_derived_1}) and (\ref{eqn:K-S_Cauchy_derived_2}) together, we get: \begin{equation}\label{eqn:K-S_Cauchy_derived_3}\frac{d}{dt}\bigg|_{t=0}(F'_tv) = -\frac{1}{2\pi i}\,\int_{\zeta\in C}(\Delta'_\omega-\zeta)^{-1}\,\bigg(A'(u_{0,\,\zeta})\bigg)\,d\zeta, \hspace{5ex} v\in C^\infty_{p,\,q}(X,\,\C),\, t\in(-\varepsilon,\,\varepsilon).\end{equation}

Let $(e'_j)_{j\in\N}$ be an orthonormal basis (w.r.t. the $L^2_\omega$-inner product) of $C^\infty_{p,\,q}(X,\,\C)$ consisting of eigenvectors of the elliptic operator $\Delta'_\omega$ corresponding respectively to the eigenvalues $(\lambda_j)_{j\in\N}$. Thus, $\Delta'_\omega e'_j = \lambda_j\,e'_j$ for every $j\in\N$. Considering the decomposition \begin{equation*}A'(u_{0,\,\zeta}) = \sum\limits_{j=0}^{+\infty}c_j(\zeta)\,e_j,  \hspace{6ex} \zeta\in D_1, \end{equation*} of the form $A'(u_{0,\,\zeta})\in C^\infty_{p,\,q}(X,\,\C)$ w.r.t. this orthonormal basis for every fixed $\zeta$, we notice that every coefficient $c_j(\zeta) = \langle\langle A'(u_{0,\,\zeta}),\,e_j\rangle\rangle_\omega$ depends holomorphically on $\zeta$ because both $u_{0,\,\zeta}=(\Delta'_\omega-\zeta)^{-1}v$ and the expression (\ref{eqn:A'u_zeta_notation}) for $A'(\cdot)$ do.

We get: \begin{equation*}(\Delta'_\omega-\zeta)^{-1}\,\bigg(A'(u_{0,\,\zeta})\bigg) =  \sum\limits_{j=0}^{+\infty}\frac{c_j(\zeta)}{\lambda_j-\zeta}\,e_j,  \hspace{6ex} \zeta\in D_1.\end{equation*} Hence,  \begin{eqnarray*}-\frac{1}{2\pi i}\,\int_{\zeta\in C}(\Delta'_\omega-\zeta)^{-1}\,\bigg(A'(u_{0,\,\zeta})\bigg)\,d\zeta = \sum\limits_{j=0}^{+\infty}\bigg(\frac{1}{2\pi i}\,\int_{\zeta\in C}\frac{c_j(\zeta)}{\zeta-\lambda_j}\,d\zeta\bigg)\,e_j = \sum\limits_{j\in J}c_j(\lambda_j)\,e_j,\end{eqnarray*} where the last equality is an application of the elementary Cauchy formula to each holomorphic function $c_j$ and $J$ is the set of indices $j$ for which $\lambda_j$ lies in the interior of the circle $C$. (Of course, the last integral above vanishes whenever $\lambda_j$ lies in the exterior of $C$.) However, we have chosen $C$ such that $0$ is the only eigenvalue of $\Delta'_\omega$ lying in the interior of $C$. Thus, $c_j(\lambda_j) = c_j(0)$ for every $j\in J$, $(e'_j)_{j\in J}$ is an $L^2_\omega$-orthonormal basis of $\ker\Delta'_\omega$ in bidegree $(p,\,q)$ and the cardinality of $J$ is the dimension of the cohomology group $H^{p,\,q}_\partial(X,\,\C)$. Therefore, we conclude that \begin{eqnarray}\label{eqn:K-S_Cauchy_derived_4}-\frac{1}{2\pi i}\,\int_{\zeta\in C}(\Delta'_\omega-\zeta)^{-1}\,\bigg(A'(u_{0,\,\zeta})\bigg)\,d\zeta = \sum\limits_{j\in J}c_j(0)\,e_j  = F'_0(A'(u_{0,\,0})),\end{eqnarray} where $F'_t:C^\infty_{p,\,q}(X,\,\C)\longrightarrow{\cal H}^{p,\,q}_{\Delta'_{\omega+t\gamma}}(X,\,\C)$ is the $L^2_\omega$-orthogonal projection onto $\ker\Delta'_\omega={\cal H}^{p,\,q}_{\Delta'_\omega}(X,\,\C)$.

Putting (\ref{eqn:K-S_Cauchy_derived_3}) and (\ref{eqn:K-S_Cauchy_derived_4}) together, we get the contention.  \hfill $\Box$

\section{Computation of the critical points of the functional $F$}\label{section:critical-points-computation_functional_F} We now pick up where we left off at the end of $\S$\ref{subsection:functional_skt-Delta}.

We fix an SKT metric $\omega$ on $X$ and we vary it along the path $\omega + t\gamma$, where $\gamma\in\ker(\partial\bar\partial)\cap C^\infty_{1,\,1}(X,\,\R)$ is an arbitrary pluriclosed real $(1,\,1)$-form. We will be considering reals $t\in(-\varepsilon,\,\varepsilon)$ with $\varepsilon>0$ small enough to ensure that $\omega + t\gamma>0$ (i.e. $\omega + t\gamma$ is a, necessarily SKT, metric on $X$). Thus, $\omega + t\gamma\in{\cal S}$ for all $t\in(-\varepsilon,\,\varepsilon)$.

For each $t\in(-\varepsilon,\,\varepsilon)$, let $\eta_{\omega,\,\gamma,\,t}$ be the minimal $L^2_{\omega + t\gamma}$-norm solution of the equation \begin{eqnarray}\label{eqn:eta_omega-gamma-t_def}d\eta = P^{(3)}_{\omega + t\gamma}(\partial\gamma),\end{eqnarray} where $P^{(3)}_{\omega + t\gamma}:C^\infty_3(X,\,\C)\longrightarrow\mbox{Im}\,d$ is the orthogonal projection w.r.t. the $L^2_{\omega + t\gamma}$-inner product. Set $\eta_{\omega,\,\gamma}=\eta_{\omega,\,\gamma,\,0}$.

Meanwhile, for each $t\in(-\varepsilon,\,\varepsilon)$, Definition \ref{Def:torsion-form} ascribes a unique torsion form $\rho_{\omega + t\gamma}\in C^\infty_2(X,\,\C)$ to $\omega + t\gamma$, identified by the property $\rho_{\omega + t\gamma}\in\mbox{Im}\,d^\star_{\omega + t\gamma}$ and the first equality below: \begin{eqnarray}\label{eqn:torsion-form_omega-t-gamma_1st-cond}\nonumber d\rho_{\omega + t\gamma} & = & P^{(3)}_{\omega + t\gamma}(\partial\omega + t\,\partial\gamma) \stackrel{(a)}{=} P^{(3)}_{\omega + t\gamma}(\partial\omega) + t\,d\eta_{\omega,\,\gamma,\,t} \\
\nonumber  & \stackrel{(b)}{=} & P^{(3)}_\omega(\partial\omega) + d(t\,\eta_{\omega,\,\gamma}) + \bigg(\frac{dP^{(3)}_{\omega + t\gamma}}{dt}_{|t=0}\bigg)(\partial\omega)\,t + {\cal O}(t^2) \\
& \stackrel{(c)}{=} & d(\rho_\omega + t\,\eta_{\omega,\,\gamma}) + \bigg(\frac{dP^{(3)}_{\omega + t\gamma}}{dt}_{|t=0}\bigg)(\partial\omega)\,t + {\cal O}(t^2),  \hspace{5ex} t\in(-\varepsilon,\,\varepsilon),\end{eqnarray} where equality (a) follows from the linearity of $P^{(3)}_{\omega + t\gamma}$, while equality (b) follows from the expansion to order two w.r.t. $t$ about $0$ of the operator $P^{(3)}_{\omega + t\gamma}$ and of the form \begin{eqnarray*}\eta_{\omega,\,\gamma,\,t} = \Delta^{-1}_{\omega +t\,\gamma}d^\star_{\omega +t\,\gamma}P^{(3)}_{\omega + t\gamma}(\partial\gamma).\end{eqnarray*} Note that both dependencies on $t$ are $C^\infty$. Equality (c) follows from (a) of (\ref{eqn:system_torsion-form}).

Recall that $F(\omega + t\gamma) = ||\rho_{\omega + t\gamma}||^2_{\omega + t\gamma}$ for every $t\in(-\varepsilon,\,\varepsilon)$. (See formula (\ref{eqn:F_energy-functional_SKT}).) The following result constitutes the first step in the computation of $\frac{d}{dt}_{|t=0} F(\omega + t\gamma) = \frac{d}{dt}_{|t=0}||\rho_{\omega + t\gamma}||^2_{\omega + t\gamma}$.

\begin{Lem}\label{Lem:1st-derivative_rho-omega+tgamma} The following equality holds: \begin{equation}\label{eqn:1st-derivative_rho-omega+tgamma}\frac{d}{dt}_{|t=0}||\rho_{\omega + t\gamma}||^2_{\omega + t\gamma}  = \frac{d}{dt}_{|t=0}||\rho_\omega + t\,\eta_{\omega,\,\gamma}||^2_{\omega + t\gamma} + 2\,||\rho_\omega||_\omega\,||\Delta_\omega^{-1}d^\star_\omega A_{\omega,\,\gamma,\,0}||_\omega,\end{equation} where \begin{equation}\label{eqn:A_omega-gamma-0_notation}A_{\omega,\,\gamma,\,0}:= \bigg(\frac{dP^{(3)}_{\omega + t\gamma}}{dt}_{|t=0}\bigg)(\partial\omega).\end{equation}

\end{Lem}

\noindent {\it Proof.} First of all we note that the formula holds if $\rho_\omega=0$. Indeed, if this is the case the right hand side is obviously zero, whereas the vanishing of the left hand side follows directly from Lemma \ref{Lem:vanishing-F-Kaehler}. Henceforth we shall assume that $\rho_{\omega}\neq 0$.

 If we put $B_{\omega,\,\gamma,\,t}:= d\rho_{\omega + t\gamma} - d(\rho_\omega + t\,\eta_{\omega,\,\gamma}) = \bigg(\frac{dP^{(3)}_{\omega + t\gamma}}{dt}_{|t=0}\bigg)(\partial\omega)\,t + {\cal O}(t^2)$ (the sum of the last two terms in (\ref{eqn:torsion-form_omega-t-gamma_1st-cond})), (\ref{eqn:torsion-form_omega-t-gamma_1st-cond}) reads: \begin{eqnarray}\label{eqn:torsion-form_omega-t-gamma_1st-cond_bis}d(\rho_\omega + t\,\eta_{\omega,\,\gamma}) = d\rho_{\omega + t\gamma} - B_{\omega,\,\gamma,\,t},  \hspace{5ex} t\in(-\varepsilon,\,\varepsilon).\end{eqnarray}

For every $t\in(-\varepsilon,\,\varepsilon)$, let $\tilde\rho_{\omega,\,\gamma,\,t}\in\mbox{Im}\,d^\star_{\omega + t\gamma}$ be the minimal $L^2_{\omega + t\gamma}$-norm solution of the equation \begin{eqnarray}\label{eqn:torsion-form_omega-t-gamma_1st-cond_bis_equation}d\rho = d\rho_{\omega + t\gamma} - B_{\omega,\,\gamma,\,t}.\end{eqnarray} This means that, for $t\in(-\varepsilon,\,\varepsilon)$, we have: $d\tilde\rho_{\omega,\,\gamma,\,t} = d\rho_{\omega + t\gamma} - B_{\omega,\,\gamma,\,t}$ and  \begin{eqnarray*}\tilde\rho_{\omega,\,\gamma,\,t} = \Delta^{-1}_{\omega + t\gamma}d^\star_{\omega + t\gamma}(d\rho_{\omega + t\gamma} - B_{\omega,\,\gamma,\,t}) = \rho_{\omega + t\gamma} - \Delta^{-1}_{\omega + t\gamma}d^\star_{\omega + t\gamma}B_{\omega,\,\gamma,\,t}.\end{eqnarray*} (To get the last equality, we used the fact that $\rho_{\omega + t\gamma}\in\mbox{Im}\,d^\star_{\omega + t\gamma}\subset\ker d^\star_{\omega + t\gamma}$.) Hence,  \begin{eqnarray}\label{eqn:torsion-form_omega-t-gamma_1st-cond_bis_proof_1}||\tilde\rho_{\omega,\,\gamma,\,t}||_{\omega + t\gamma}\geq||\rho_{\omega + t\gamma}||_{\omega + t\gamma} - ||\Delta^{-1}_{\omega + t\gamma}d^\star_{\omega + t\gamma}B_{\omega,\,\gamma,\,t}||_{\omega + t\gamma},  \hspace{5ex} t\in(-\varepsilon,\,\varepsilon).\end{eqnarray}

On the other hand, the $L^2_{\omega + t\gamma}$-norm minimality of $\tilde\rho_{\omega,\,\gamma,\,t}\in\mbox{Im}\,d^\star_{\omega + t\gamma}$ among the solutions of equation (\ref{eqn:torsion-form_omega-t-gamma_1st-cond_bis_equation}), one of which is $\rho_\omega + t\,\eta_{\omega,\,\gamma}$ (cf. (\ref{eqn:torsion-form_omega-t-gamma_1st-cond_bis})), implies that $||\rho_\omega + t\,\eta_{\omega,\,\gamma}||_{\omega + t\gamma}\geq||\tilde\rho_{\omega,\,\gamma,\,t}||_{\omega + t\gamma}$ for each $t\in(-\varepsilon,\,\varepsilon)$. Together with (\ref{eqn:torsion-form_omega-t-gamma_1st-cond_bis_proof_1}), this yileds: \begin{eqnarray}\label{eqn:torsion-form_omega-t-gamma_1st-cond_bis_proof_2}||\rho_\omega + t\,\eta_{\omega,\,\gamma}||_{\omega + t\gamma} - ||\rho_{\omega + t\gamma}||_{\omega + t\gamma} + ||\Delta^{-1}_{\omega + t\gamma}d^\star_{\omega + t\gamma}B_{\omega,\,\gamma,\,t}||_{\omega + t\gamma}\geq 0,  \hspace{5ex} t\in(-\varepsilon,\,\varepsilon).\end{eqnarray} Note that the l.h.s. term of (\ref{eqn:torsion-form_omega-t-gamma_1st-cond_bis_proof_2}) vanishes at $t=0$, so, as a function of $t\in(-\varepsilon,\,\varepsilon)$, it achieves its minimum at $t=0$. In particular, its derivative w.r.t. $t$ vanishes at $t=0$. Here, we make use of the assumption $\rho_\omega\neq 0$ (otherwise, $\omega$ is K\"ahler -- see Lemma \ref{Lem:Kaehler_torsion-vanishing} -- and there is nothing to prove), so the norms are smooth and the differentiation is justified. This yields: \begin{eqnarray}\label{eqn:torsion-form_omega-t-gamma_1st-cond_bis_proof_3}\frac{d}{dt}_{|t=0}||\rho_\omega + t\,\eta_{\omega,\,\gamma}||_{\omega + t\gamma} = \frac{d}{dt}_{|t=0}||\rho_{\omega + t\gamma}||_{\omega + t\gamma} - \frac{d}{dt}_{|t=0}||\Delta^{-1}_{\omega + t\gamma}d^\star_{\omega + t\gamma}B_{\omega,\,\gamma,\,t}||_{\omega + t\gamma}.\end{eqnarray}

Meanwhile, $\frac{d}{dt}_{|t=0}||\rho_{\omega + t\gamma}||_{\omega + t\gamma}^2 = 2\,||\rho_\omega||_\omega\,\frac{d}{dt}_{|t=0}\,||\rho_{\omega + t\gamma}||_{\omega + t\gamma}$. Thus, using (\ref{eqn:torsion-form_omega-t-gamma_1st-cond_bis_proof_3}), we get:

\begin{eqnarray}\label{eqn:torsion-form_omega-t-gamma_1st-cond_bis_proof_4}\nonumber\frac{d}{dt}_{|t=0}||\rho_{\omega + t\gamma}||_{\omega + t\gamma}^2 & = & 2\,||\rho_\omega||_\omega\, \frac{d}{dt}_{|t=0}||\rho_\omega + t\,\eta_{\omega,\,\gamma}||_{\omega + t\gamma} + 2\,||\rho_\omega||_\omega\,\frac{d}{dt}_{|t=0}\bigg(||\Delta^{-1}_{\omega + t\gamma}d^\star_{\omega + t\gamma}B_{\omega,\,\gamma,\,t}||_{\omega + t\gamma}^2\bigg)^{\frac{1}{2}} \\
  & = & \frac{d}{dt}_{|t=0}||\rho_\omega + t\,\eta_{\omega,\,\gamma}||_{\omega + t\gamma}^2 + 2\,||\rho_\omega||_\omega\,M(0),\end{eqnarray} where we put $$M(t):= \frac{\frac{d}{dt}\,||\Delta^{-1}_{\omega + t\gamma}d^\star_{\omega + t\gamma}B_{\omega,\,\gamma,\,t}||_{\omega + t\gamma}^2}{2\,||\Delta^{-1}_{\omega + t\gamma}d^\star_{\omega + t\gamma}B_{\omega,\,\gamma,\,t}||_{\omega + t\gamma}},  \hspace{5ex} t\in(-\varepsilon,\,\varepsilon).$$

Now, letting $A_{\omega,\,\gamma,\,t} = \frac{1}{t}\, B_{\omega,\,\gamma,\,t}$, we see that $A_{\omega,\,\gamma,\,0}$ is given by the expression (\ref{eqn:A_omega-gamma-0_notation}) and we have $$A_{\omega,\,\gamma,\,t} = \bigg(\frac{dP^{(3)}_{\omega + t\gamma}}{dt}_{|t=0}\bigg)(\partial\omega) + {\cal O}(t),  \hspace{5ex} t\in(-\varepsilon,\,\varepsilon).$$ Thus, writing $B_{\omega,\,\gamma,\,t} = t\,A_{\omega,\,\gamma,\,t}$ in the expression for $M(t)$, we get for $t\in(0,\,\varepsilon)$: \begin{eqnarray*}M(t)= \frac{\frac{d}{dt}\,(t^2\,||\Delta^{-1}_{\omega + t\gamma}d^\star_{\omega + t\gamma}A_{\omega,\,\gamma,\,t}||_{\omega + t\gamma}^2)}{2t\,||\Delta^{-1}_{\omega + t\gamma}d^\star_{\omega + t\gamma}A_{\omega,\,\gamma,\,t}||_{\omega + t\gamma}} = ||\Delta^{-1}_{\omega + t\gamma}d^\star_{\omega + t\gamma}A_{\omega,\,\gamma,\,t}||_{\omega + t\gamma} + \frac{t}{2}\,\frac{\frac{d}{dt}\,||\Delta^{-1}_{\omega + t\gamma}d^\star_{\omega + t\gamma}A_{\omega,\,\gamma,\,t}||_{\omega + t\gamma}^2}{||\Delta^{-1}_{\omega + t\gamma}d^\star_{\omega + t\gamma}A_{\omega,\,\gamma,\,t}||_{\omega + t\gamma}}.\end{eqnarray*} Hence, $M(0)= ||\Delta^{-1}_{\omega}d^\star_{\omega }A_{\omega,\,\gamma,\,0}||_{\omega}$. Together with (\ref{eqn:torsion-form_omega-t-gamma_1st-cond_bis_proof_4}), this proves (\ref{eqn:1st-derivative_rho-omega+tgamma}). \hfill $\Box$

\vspace{3ex}

Since $F(\omega + t\gamma) = ||\rho_{\omega + t\gamma}||^2_{\omega + t\gamma}$ for every $t\in(-\varepsilon,\,\varepsilon)$ (cf. formula (\ref{eqn:F_energy-functional_SKT})), the above Lemma \ref{Lem:1st-derivative_rho-omega+tgamma} reduces the computation of $\frac{d}{dt}_{|t=0} F(\omega + t\gamma) = \frac{d}{dt}_{|t=0}||\rho_{\omega + t\gamma}||^2_{\omega + t\gamma}$ to the computation of $\frac{d}{dt}_{|t=0}||\rho_\omega + t\,\eta_{\omega,\,\gamma}||^2_{\omega + t\gamma}$. This last derivative at $t=0$ is computed in the following

\begin{Lem}\label{Lem:1st-variation_Ft-tilde} Let $X$ be a compact complex manifold with $\mbox{dim}_\C X=n$. Suppose there exists an SKT metric $\omega$ on $X$.

  For any real $(1,\,1)$-form $\gamma\in\ker(\partial\bar\partial)$ and for a real $t$ varying in a neighbourhood $(-\varepsilon,\,\varepsilon)$ of $0$ so small that $\omega + t\gamma>0$, we have \begin{eqnarray}\label{eqn:1st-variation_Ft-tilde}\nonumber\frac{d}{dt}_{|t=0}||\rho_\omega + t\,\eta_{\omega,\,\gamma}||^2_{\omega + t\gamma} & = & \langle\langle\eta_{\omega,\,\gamma}^{2,\,0},\,\rho_\omega^{2,\,0}\rangle\rangle_\omega + \langle\langle\rho_\omega^{2,\,0},\,\eta_{\omega,\,\gamma}^{2,\,0} + [\Lambda_\omega,\,\gamma\wedge\cdot]\,\rho_\omega^{2,\,0}\rangle\rangle_\omega \\
\nonumber & + & \langle\langle\eta_{\omega,\,\gamma}^{0,\,2},\,\rho_\omega^{0,\,2}\rangle\rangle_\omega + \langle\langle\rho_\omega^{0,\,2},\,\eta_{\omega,\,\gamma}^{0,\,2} + [\Lambda_\omega,\,\gamma\wedge\cdot]\,\rho_\omega^{0,\,2}\rangle\rangle_\omega \\    
     & + & \langle\langle\eta_{\omega,\,\gamma}^{1,\,1},\,\rho_\omega^{1,\,1}\rangle\rangle_\omega + \langle\langle\rho_\omega^{1,\,1},\,\eta_{\omega,\,\gamma}^{1,\,1} + [\Lambda_\omega,\,\gamma\wedge\cdot]\,\rho_\omega^{1,\,1}\rangle\rangle_\omega,\end{eqnarray} where $\rho_\omega = \rho_\omega^{2,\,0} + \rho_\omega^{1,\,1} + \rho_\omega^{0,\,2}$ and $\eta_{\omega,\,\gamma} = \eta_{\omega,\,\gamma}^{2,\,0} + \eta_{\omega,\,\gamma}^{1,\,1} + \eta_{\omega,\,\gamma}^{0,\,2}$ are the decompositions of $\rho_\omega$ and $\eta_{\omega,\,\gamma}$ into pure-type forms. 

\end{Lem}

\noindent {\it Proof.} For every $(-\varepsilon,\,\varepsilon)$, let $\alpha_t:= \rho_\omega + t\,\eta_{\omega,\,\gamma}$. In particular, $\alpha_0 = \rho_\omega$.

Since pure-type forms of different types are mutually orthogonal, we have for all $t$ near $0$: \begin{eqnarray}\label{eqn:1st-variation_Ft-tilde_proof}\nonumber||\alpha_t||^2_{\omega + t\gamma} & = & \int\limits_X|\alpha_t^{2,\,0}|^2_{\omega + t\gamma}\,dV_{\omega + t\gamma} + \int\limits_X|\alpha_t^{0,\,2}|^2_{\omega + t\gamma}\,dV_{\omega + t\gamma} + \int\limits_X|\alpha_t^{1,\,1} |^2_{\omega + t\gamma}\,dV_{\omega + t\gamma} \\
  & = & \int\limits_X\alpha_t^{2,\,0}\wedge\star_{\omega + t\gamma}\overline{\alpha_t^{2,\,0}} + \int\limits_X\alpha_t^{0,\,2}\wedge\star_{\omega + t\gamma}\overline{\alpha_t^{0,\,2}} + \int\limits_X\alpha_t^{1,\,1}\wedge\star_{\omega + t\gamma}\overline{\alpha_t^{1,\,1}}.\end{eqnarray}

Now, for all $(p,\,q)\in\{(2,\,0),\,(1,\,1),\,(0,\,2)\}$, we have: $\frac{d\alpha_t^{p,\,q}}{dt}_{|t=0} = \eta_{\omega,\,\gamma}^{p,\,q}$ and, thanks to Lemma \ref{Lem:1st-variation_Hodge-star},   \begin{eqnarray*}\frac{d(\star_{\omega + t\gamma}\overline{\alpha_t^{p,\,q}})}{dt}_{\bigg|t=0} = \star_\omega\bigg(\overline{\eta_{\omega,\,\gamma}^{p,\,q}} + [\Lambda_\omega,\,\gamma\wedge\cdot]\,\overline{\rho_\omega^{p,\,q}}\bigg).\end{eqnarray*}

After differentiating (\ref{eqn:1st-variation_Ft-tilde_proof}) with respect to $t$ and then evaluating the result at $t=0$, formula (\ref{eqn:1st-variation_Ft-tilde}) follows at once from the above equalities thanks to the equality $u\wedge\star_\omega\bar{v} = \langle u,\,v\rangle_\omega\,dV_\omega$ holding for all forms $u$ and $v$ of the same degree and defining the Hodge star operator.  \hfill $\Box$

\vspace{3ex}

Putting Lemmas \ref{Lem:1st-derivative_rho-omega+tgamma} and \ref{Lem:1st-variation_Ft-tilde} together, we get the following

\begin{The}\label{The:1st-variation_Ft-tilde_final} Let $X$ be a compact complex manifold with $\mbox{dim}_\C X=n$. Suppose there exists an SKT metric $\omega$ on $X$.

  For any real $(1,\,1)$-form $\gamma\in\ker(\partial\bar\partial)$ and for a real $t$ varying in a neighbourhood $(-\varepsilon,\,\varepsilon)$ of $0$ so small that $\omega + t\gamma>0$, we have: \begin{eqnarray}\label{eqn:1st-variation_Ft-tilde_final}\nonumber (d_\omega F)(\gamma) = \frac{d}{dt}_{|t=0} F(\omega + t\gamma)  & = & \langle\langle\eta_{\omega,\,\gamma}^{2,\,0},\,\rho_\omega^{2,\,0}\rangle\rangle_\omega + \langle\langle\rho_\omega^{2,\,0},\,\eta_{\omega,\,\gamma}^{2,\,0} + [\Lambda_\omega,\,\gamma\wedge\cdot]\,\rho_\omega^{2,\,0}\rangle\rangle_\omega \\
\nonumber & + & \langle\langle\eta_{\omega,\,\gamma}^{0,\,2},\,\rho_\omega^{0,\,2}\rangle\rangle_\omega + \langle\langle\rho_\omega^{0,\,2},\,\eta_{\omega,\,\gamma}^{0,\,2} + [\Lambda_\omega,\,\gamma\wedge\cdot]\,\rho_\omega^{0,\,2}\rangle\rangle_\omega \\    
 \nonumber    & + & \langle\langle\eta_{\omega,\,\gamma}^{1,\,1},\,\rho_\omega^{1,\,1}\rangle\rangle_\omega + \langle\langle\rho_\omega^{1,\,1},\,\eta_{\omega,\,\gamma}^{1,\,1} + [\Lambda_\omega,\,\gamma\wedge\cdot]\,\rho_\omega^{1,\,1}\rangle\rangle_\omega\\
 & + & 2\,||\rho_\omega||_\omega\,||\Delta_\omega^{-1}d^\star_\omega A_{\omega,\,\gamma,\,0}||_\omega, \end{eqnarray} where $A_{\omega,\,\gamma,\,0}$ is defined by formula (\ref{eqn:A_omega-gamma-0_notation}).

In particular, if we choose $\gamma = \omega$, the above formula specialises to \begin{eqnarray}\label{eqn:1st-variation_Ft-tilde_final_special}(d_\omega F)(\gamma) = \frac{d}{dt}_{|t=0} F(\omega + t\gamma)  & = & n\,||\rho_\omega||^2_\omega.\end{eqnarray}

\end{The}

\noindent {\it Proof.} Only formula (\ref{eqn:1st-variation_Ft-tilde_final_special}) still needs a proof. When $\gamma = \omega$, we have $\eta_{\omega,\,\gamma} = \rho_\omega$ (see formulae (\ref{eqn:eta_omega-gamma-t_def}) and (\ref{eqn:system_torsion-form})) and $[\Lambda_\omega,\,\gamma\wedge\cdot] = [\Lambda_\omega,\,\omega\wedge\cdot] = (n-p-q)\,\mbox{Id}$ on $(p,\,q)$-forms for all bidegrees $(p,\,q)$ (see e.g. [Voi02, Lemma 6.19] for the last, standard, equality). Thus, the quantities on the first three lines of the right-hand side of (\ref{eqn:1st-variation_Ft-tilde_final}) are equal to $n\,||\rho_\omega^{2,\,0}||^2_\omega$, $n\,||\rho_\omega^{0,\,2}||^2_\omega$ and respectively $n\,||\rho_\omega^{1,\,1}||^2_\omega$ when $\gamma = \omega$.

On the other hand, using (\ref{eqn:A_omega-gamma-0_notation}), we get: \begin{eqnarray*}A_{\omega,\,\omega,\,0} = \bigg(\frac{dP^{(3)}_{(1+t)\,\omega}}{dt}_{|t=0}\bigg)(\partial\omega) = 0 \end{eqnarray*} since $P^{(3)}_{(1+t)\,\omega}(\partial\omega) =  P^{(3)}_\omega(\partial\omega)$ for every $t$ close enough to $0$. (See (\ref{eqn:P^3_del-omega_lambda}) for this last fact.)

Putting these pieces of information together, we get (\ref{eqn:1st-variation_Ft-tilde_final_special}).  \hfill $\Box$

\vspace{3ex}

Note that formula (\ref{eqn:1st-variation_Ft-tilde_final_special}) coincides with formula (\ref{eqn:rescaling_F-critical_1+t}) of Lemma \ref{Lem:rescaling_F-critical}. Another conclusion of this study is the following result that was already proved in a different way in (\ref{eqn:critical-points_F-Kaehler}) of Lemma \ref{Lem:rescaling_F-critical}.

\begin{Cor}\label{Cor:crit-points_SKT-primitive-torsion} Let $X$ be a compact complex SKT manifold.

Then, for every SKT metric $\omega$ on $X$, $\omega$ is a {\bf critical point} for the functional $F:{\cal S}\longrightarrow[0,\,+\infty)$ if and only if $\omega$ is a {\bf K\"ahler} metric.

\end{Cor}  

\noindent {\it Proof.} If an SKT metric $\omega$ is a critical point for $F$, then $(d_\omega F)(\gamma) = \frac{d}{dt}_{|t=0}F(\omega + t\,\gamma) = 0$ for every real $(1,\,1)$-form $\gamma$ lying in $\ker(\partial\bar\partial)$. In particular, choosing $\gamma:=\omega$ (which is allowed since $\partial\bar\partial\omega=0$), the vanishing of $\frac{d}{dt}_{|t=0}F(\omega + t\,\omega)$ implies the vanishing of $\rho_\omega$ thanks to (\ref{eqn:1st-variation_Ft-tilde_final_special}). On the other hand, the equality $\rho_\omega=0$ implies that $\omega$ is K\"ahler by Lemma \ref{Lem:Kaehler_torsion-vanishing}.

Conversely, if $\omega$ is a K\"ahler metric, $\omega$ is a minimiser, hence a critical point, for the functional $F$ by Lemma \ref{Lem:vanishing-F-Kaehler}. 

One can also argue that $\rho_\omega=0$ if $\omega$ is K\"ahler (again by Lemma \ref{Lem:Kaehler_torsion-vanishing}), so (\ref{eqn:1st-variation_Ft-tilde_final}) of Theorem \ref{The:1st-variation_Ft-tilde_final} implies that $\frac{d}{dt}_{|t=0}F(\omega + t\,\gamma) = 0$ for every real $(1,\,1)$-form $\gamma$ lying in $\ker(\partial\bar\partial)$. This means that $\omega$ is a critical point for $F$.  \hfill $\Box$

\section{Computation of the critical points of the functional $G$}\label{section:critical-points-computation_functional_G} We now pick up where we left off at the end of $\S$\ref{subsection:functional_balanced}.

We fix a balanced metric $\omega$ on $X$ (so $\omega_{n-1}\in{\cal B}$) and we vary it along the path $\omega + t\Omega$, where $\Omega\in C^\infty_{n-1,\,n-1}(X,\,\R)$ is an arbitrary real $(n-1,\,n-1)$-form such that $\bar\partial\Omega = 0$. We will be considering reals $t\in(-\varepsilon,\,\varepsilon)$ with $\varepsilon>0$ small enough to ensure that $\omega_{n-1} + t\Omega>0$ (i.e. $\omega_{n-1} + t\Omega$ is a, necessarily balanced, metric on $X$). Thus, $\omega_{n-1} + t\Omega\in{\cal B}$ for all $t\in(-\varepsilon,\,\varepsilon)$. By [Mic83], there exists a unique $C^\infty$ positive definite $(1,\,1)$-form $\gamma_{\omega,\,\Omega,\,t}$ such that \begin{equation}\label{eqn:root_balanced}(\gamma_{\omega,\,\Omega,\,t})_{n-1} = \omega_{n-1} + t\Omega,  \hspace{3ex} t\in(-\varepsilon,\,\varepsilon).\end{equation} Note that $\gamma_{\omega,\,\Omega,\,0} = \omega$. 

For each $t\in(-\varepsilon,\,\varepsilon)$, let $\eta_{\omega,\,\Omega,\,t}$ be the minimal $L^2_{\gamma_{\omega,\,\Omega,\,t}}$-norm solution of the equation \begin{eqnarray}\label{eqn:eta_omega-Omega-t_def}\bar\partial\eta_{\omega,\,\Omega,\,t} = P^{(n-1,\,n-1)}_{\omega_{n-1} + t\Omega}(\Omega),\end{eqnarray} where $P^{(n-1,\,n-1)}_{\omega_{n-1} + t\Omega}:C^\infty_{n-1,\,n-1}(X,\,\C)\longrightarrow\mbox{Im}\,\bar\partial$ is the orthogonal projection w.r.t. the $L^2_{\gamma_{\omega,\,\Omega,\,t}}$-inner product. Set $\eta_{\omega,\,\Omega}=\eta_{\omega,\,\Omega,\,0}$.

Meanwhile, for each $t\in(-\varepsilon,\,\varepsilon)$, Definition \ref{Def:balanced_torsion-form} ascribes a unique torsion $(n-1,\,n-2)$-form $\Gamma_{\omega_{n-1} + t\Omega}$ to $\omega_{n-1} + t\Omega$, identified by the property $\Gamma_{\omega_{n-1} + t\Omega}\in\mbox{Im}\,\bar\partial^\star_{\gamma_{\omega,\,\Omega,\,t}}$ and the first equality below: \begin{eqnarray}\label{eqn:torsion-form_omega-Omega-t_1st-cond}\nonumber \bar\partial\Gamma_{\omega_{n-1} + t\Omega} & = & P^{(n-1,\,n-1)}_{\omega_{n-1} + t\Omega}(\omega_{n-1} + t\Omega) \stackrel{(a)}{=} P^{(n-1,\,n-1)}_{\omega_{n-1} + t\Omega}(\omega_{n-1}) + t\,\bar\partial\eta_{\omega,\,\Omega,\,t} \\
\nonumber  & \stackrel{(b)}{=} & P^{(n-1,\,n-1)}_{\omega_{n-1}}(\omega_{n-1}) + t\,\bar\partial\eta_{\omega,\,\Omega} + \bigg(\frac{dP^{(n-1,\,n-1)}_{\omega_{n-1} + t\Omega}}{dt}_{|t=0}\bigg)(\omega_{n-1})\,t + {\cal O}(t^2) \\
& \stackrel{(c)}{=} & \bar\partial(\Gamma_{\omega_{n-1}} + t\,\eta_{\omega,\,\Omega}) + \bigg(\frac{dP^{(n-1,\,n-1)}_{\omega_{n-1} + t\Omega}}{dt}_{|t=0}\bigg)(\omega_{n-1})\,t + {\cal O}(t^2),  \hspace{5ex} t\in(-\varepsilon,\,\varepsilon),\end{eqnarray} where equality (a) follows from the linearity of $P^{(n-1,\,n-1)}_{\omega_{n-1} + t\Omega}$, while equality (b) follows from the expansion to order two w.r.t. $t$ about $0$ of the operator $P^{(3)}_{\omega + t\gamma}$ and of the form \begin{eqnarray*}\eta_{\omega,\,\Omega,\,t} = \Delta^{''-1}_{\gamma_{\omega,\,\Omega,\,t}}\bar\partial^\star_{\gamma_{\omega,\,\Omega,\,t}}P^{(n-1,\,n-1)}_{\omega_{n-1} + t\Omega}(\Omega).\end{eqnarray*} Note that both dependencies on $t$ are $C^\infty$. Equality (c) follows from (a) of (\ref{eqn:system_balanced_torsion-form}).

Recall that $G(\omega_{n-1} + t\Omega) = ||\Gamma_{\omega_{n-1} + t\Omega}||^2_{\gamma_{\omega,\,\Omega,\,t}}$ for every $t\in(-\varepsilon,\,\varepsilon)$. (See formula (\ref{eqn:G_energy-functional_balanced}).) The following result constitutes the first step in the computation of $\frac{d}{dt}_{|t=0} G(\omega_{n-1} + t\Omega) = \frac{d}{dt}_{|t=0}||\Gamma_{\omega_{n-1} + t\Omega}||^2_{\gamma_{\omega,\,\Omega,\,t}}$.

\begin{Lem}\label{Lem:1st-derivative_Gamma_omega_n-1+tOmega} The following equality holds: \begin{equation}\label{eqn:1st-derivative_Gamma_omega_n-1+tOmega}\frac{d}{dt}_{|t=0}||\Gamma_{\omega_{n-1} + t\Omega}||^2_{\gamma_{\omega,\,\Omega,\,t}}  = \frac{d}{dt}_{|t=0}||\Gamma_{\omega_{n-1}} + t\,\eta_{\omega,\,\Omega}||^2_{\gamma_{\omega,\,\Omega,\,t}} + 2\,||\Gamma_{\omega_{n-1}}||_\omega\,||\Delta_\omega^{''-1}\bar\partial^\star_\omega A_{\omega,\,\Omega,\,0}||_\omega,\end{equation} where \begin{equation}\label{eqn:A_omega-Omega-0_notation}A_{\omega,\,\Omega,\,0}:= \bigg(\frac{dP^{(n-1,\,n-1)}_{\omega_{n-1} + t\Omega}}{dt}_{|t=0}\bigg)(\omega_{n-1}).\end{equation}

\end{Lem}

\noindent {\it Proof.} Again if $\Gamma_{\omega_{n-1}}=0$ the formula holds. Thus we may assume from now on that $\Gamma_{\omega_{n-1}}\neq 0$.If we put $B_{\omega,\,\Omega,\,t}:= \bar\partial\Gamma_{\omega_{n-1} + t\Omega} - \bar\partial(\Gamma_{\omega_{n-1}} + t\,\eta_{\omega,\,\Omega}) = \bigg(\frac{dP^{(n-1,\,n-1)}_{\omega_{n-1} + t\Omega}}{dt}_{|t=0}\bigg)(\omega_{n-1})\,t + {\cal O}(t^2)$ (the sum of the last two terms in (\ref{eqn:torsion-form_omega-Omega-t_1st-cond})), (\ref{eqn:torsion-form_omega-Omega-t_1st-cond}) reads: \begin{eqnarray}\label{eqn:torsion-form_omega-Omega-t_1st-cond_bis}\bar\partial(\Gamma_{\omega_{n-1}} + t\,\eta_{\omega,\,\Omega}) = \bar\partial\Gamma_{\omega_{n-1} + t\Omega} - B_{\omega,\,\Omega,\,t},  \hspace{5ex} t\in(-\varepsilon,\,\varepsilon).\end{eqnarray}

For every $t\in(-\varepsilon,\,\varepsilon)$, let $\tilde\Gamma_{\omega,\,\Omega,\,t}\in\mbox{Im}\,\bar\partial^\star_{\gamma_{\omega,\,\Omega,\,t}}$ be the minimal $L^2_{\gamma_{\omega,\,\Omega,\,t}}$-norm solution of the equation \begin{eqnarray}\label{eqn:torsion-form_omega-Omega-t_1st-cond_bis_equation}\bar\partial\Gamma = \bar\partial\Gamma_{\omega_{n-1} + t\Omega} - B_{\omega,\,\Omega,\,t}.\end{eqnarray} This means that, for $t\in(-\varepsilon,\,\varepsilon)$, we have: $\bar\partial\tilde\Gamma_{\omega,\,\Omega,\,t} = \bar\partial\Gamma_{\omega_{n-1} + t\Omega} - B_{\omega,\,\Omega,\,t}$ and  \begin{eqnarray*}\tilde\Gamma_{\omega,\,\Omega,\,t} = \Delta^{''-1}_{\gamma_{\omega,\,\Omega,\,t}}\bar\partial^\star_{\gamma_{\omega,\,\Omega,\,t}}(\bar\partial\Gamma_{\omega_{n-1} + t\Omega} - B_{\omega,\,\Omega,\,t}) = \Gamma_{\omega_{n-1} + t\Omega} - \Delta^{''-1}_{\gamma_{\omega,\,\Omega,\,t}}\bar\partial^\star_{\gamma_{\omega,\,\Omega,\,t}}B_{\omega,\,\Omega,\,t}.\end{eqnarray*} (To get the last equality, we used the fact that $\Gamma_{\omega_{n-1} + t\Omega}\in\mbox{Im}\,\bar\partial^\star_{\gamma_{\omega,\,\Omega,\,t}}\subset\ker\bar\partial^\star_{\gamma_{\omega,\,\Omega,\,t}}$.) Hence,  \begin{eqnarray}\label{eqn:torsion-form_omega-Omega-t_1st-cond_bis_proof_1}||\tilde\Gamma_{\omega,\,\Omega,\,t}||_{\gamma_{\omega,\,\Omega,\,t}}\geq||\Gamma_{\omega_{n-1} + t\Omega}||_{\gamma_{\omega,\,\Omega,\,t}} - ||\Delta^{''-1}_{\gamma_{\omega,\,\Omega,\,t}}\bar\partial^\star_{\gamma_{\omega,\,\Omega,\,t}}B_{\omega,\,\Omega,\,t}||_{\gamma_{\omega,\,\Omega,\,t}},  \hspace{5ex} t\in(-\varepsilon,\,\varepsilon).\end{eqnarray}

On the other hand, the $L^2_{\gamma_{\omega,\,\Omega,\,t}}$-norm minimality of $\tilde\Gamma_{\omega,\,\Omega,\,t}\in\mbox{Im}\,\bar\partial^\star_{\gamma_{\omega,\,\Omega,\,t}}$ among the solutions of equation (\ref{eqn:torsion-form_omega-Omega-t_1st-cond_bis_equation}), one of which is $\Gamma_{\omega_{n-1}} + t\,\eta_{\omega,\,\Omega}$ (cf. (\ref{eqn:torsion-form_omega-Omega-t_1st-cond_bis})), implies that $||\Gamma_{\omega_{n-1}} + t\,\eta_{\omega,\,\Omega}||_{\gamma_{\omega,\,\Omega,\,t}}\geq||\tilde\Gamma_{\omega,\,\Omega,\,t}||_{\gamma_{\omega,\,\Omega,\,t}}$ for each $t\in(-\varepsilon,\,\varepsilon)$. Together with (\ref{eqn:torsion-form_omega-Omega-t_1st-cond_bis_proof_1}), this yileds: \begin{eqnarray}\label{eqn:torsion-form_omega-Omega-t_1st-cond_bis_proof_2}||\Gamma_{\omega_{n-1}} + t\,\eta_{\omega,\,\Omega}||_{\gamma_{\omega,\,\Omega,\,t}} - ||\Gamma_{\omega_{n-1} + t\Omega}||_{\gamma_{\omega,\,\Omega,\,t}} + ||\Delta^{''-1}_{\gamma_{\omega,\,\Omega,\,t}}\bar\partial^\star_{\gamma_{\omega,\,\Omega,\,t}}B_{\omega,\,\Omega,\,t}||_{\gamma_{\omega,\,\Omega,\,t}}\geq 0,  \hspace{5ex} t\in(-\varepsilon,\,\varepsilon).\end{eqnarray} Note that the l.h.s. term of (\ref{eqn:torsion-form_omega-Omega-t_1st-cond_bis_proof_2}) vanishes at $t=0$, so, as a function of $t\in(-\varepsilon,\,\varepsilon)$, it achieves its minimum at $t=0$. In particular, its derivative w.r.t. $t$ vanishes at $t=0$. Once again, the assumption $\Gamma_{\omega_{n-1}}\neq0$ (otherwise, $\omega$ is K\"ahler and there is nothing to prove -- see Lemma \ref{Lem:Kaehler_balanced_torsion-vanishing}) we have made justifies the differentiability. This yields: \begin{eqnarray}\label{eqn:torsion-form_omega-Omega-t_1st-cond_bis_proof_3}\frac{d}{dt}_{|t=0}||\Gamma_{\omega_{n-1}} + t\,\eta_{\omega,\,\Omega}||_{\gamma_{\omega,\,\Omega,\,t}} = \frac{d}{dt}_{|t=0}||\Gamma_{\omega_{n-1} + t\Omega}||_{\gamma_{\omega,\,\Omega,\,t}} - \frac{d}{dt}_{|t=0}||\Delta^{''-1}_{\gamma_{\omega,\,\Omega,\,t}}\bar\partial^\star_{\gamma_{\omega,\,\Omega,\,t}}B_{\omega,\,\Omega,\,t}||_{\gamma_{\omega,\,\Omega,\,t}}.\end{eqnarray}

Meanwhile, $\frac{d}{dt}_{|t=0}||\Gamma_{\omega_{n-1} + t\Omega}||_{\gamma_{\omega,\,\Omega,\,t}}^2 = 2\,||\Gamma_{\omega_{n-1}}||_\omega\,\frac{d}{dt}_{|t=0}\,||\Gamma_{\omega_{n-1} + t\Omega}||_{\gamma_{\omega,\,\Omega,\,t}}$. Thus, using (\ref{eqn:torsion-form_omega-Omega-t_1st-cond_bis_proof_3}), we get:

\begin{eqnarray}\label{eqn:torsion-form_omega-Omega-t_1st-cond_bis_proof_4}\nonumber\frac{d}{dt}_{|t=0}||\Gamma_{\omega_{n-1} + t\Omega}||_{\gamma_{\omega,\,\Omega,\,t}}^2 & = & 2\,||\Gamma_{\omega_{n-1}}||_\omega\,\frac{d}{dt}_{|t=0}||\Gamma_{\omega_{n-1}} + t\,\eta_{\omega,\,\Omega}||_{\gamma_{\omega,\,\Omega,\,t}} \\
\nonumber  &  & + 2\,||\Gamma_{\omega_{n-1}}||_\omega\,\frac{d}{dt}_{|t=0}\bigg(||\Delta^{''-1}_{\gamma_{\omega,\,\Omega,\,t}}\bar\partial^\star_{\gamma_{\omega,\,\Omega,\,t}}B_{\omega,\,\Omega,\,t}||_{\gamma_{\omega,\,\Omega,\,t}}^2\bigg)^{\frac{1}{2}} \\
& = & \frac{d}{dt}_{|t=0}||\Gamma_{\omega_{n-1}} + t\,\eta_{\omega,\,\Omega}||_{\gamma_{\omega,\,\Omega,\,t}}^2 + 2\,||\Gamma_{\omega_{n-1}}||_\omega\,M(0),\end{eqnarray} where we put $$M(t):= \frac{\frac{d}{dt}\,||\Delta^{''-1}_{\gamma_{\omega,\,\Omega,\,t}}\bar\partial^\star_{\gamma_{\omega,\,\Omega,\,t}}B_{\omega,\,\Omega,\,t}||_{\gamma_{\omega,\,\Omega,\,t}}^2}{2\,||\Delta^{''-1}_{\gamma_{\omega,\,\Omega,\,t}}\bar\partial^\star_{\gamma_{\omega,\,\Omega,\,t}}B_{\omega,\,\Omega,\,t}||_{\gamma_{\omega,\,\Omega,\,t}}},  \hspace{5ex} t\in(-\varepsilon,\,\varepsilon).$$

Now, letting $A_{\omega,\,\Omega,\,t} = \frac{1}{t}\, B_{\omega,\,\Omega,\,t}$, we see that $A_{\omega,\,\Omega,\,0}$ is given by the expression (\ref{eqn:A_omega-Omega-0_notation}) and we have $$A_{\omega,\,\gamma,\,t} = \bigg(\frac{dP^{(n-1,\,n-1)}_{\omega_{n-1} + t\Omega}}{dt}_{|t=0}\bigg)(\omega_{n-1}) + {\cal O}(t),  \hspace{5ex} t\in(-\varepsilon,\,\varepsilon).$$ Thus, writing $B_{\omega,\,\Omega,\,t} = t\,A_{\omega,\,\Omega,\,t}$ in the expression for $M(t)$, we get for $t\in(0,\,\varepsilon)$: \begin{eqnarray*}M(t)= \frac{\frac{d}{dt}\,(t^2\,||\Delta^{''-1}_{\gamma_{\omega,\,\Omega,\,t}}\bar\partial^\star_{\gamma_{\omega,\,\Omega,\,t}}A_{\omega,\,\Omega,\,t}||_{\gamma_{\omega,\,\Omega,\,t}}^2)}{2t\,||\Delta^{''-1}_{\omega + t\gamma}\bar\partial^\star_{\gamma_{\omega,\,\Omega,\,t}}A_{\omega,\,\Omega,\,t}||_{\gamma_{\omega,\,\Omega,\,t}}} & = & ||\Delta^{''-1}_{\gamma_{\omega,\,\Omega,\,t}}\bar\partial^\star_{\gamma_{\omega,\,\Omega,\,t}}A_{\omega,\,\Omega,\,t}||_{\gamma_{\omega,\,\Omega,\,t}} \\
  & + & \frac{t}{2}\,\frac{\frac{d}{dt}\,||\Delta^{''-1}_{\gamma_{\omega,\,\Omega,\,t}}\bar\partial^\star_{\gamma_{\omega,\,\Omega,\,t}}A_{\omega,\,\Omega,\,t}||_{\gamma_{\omega,\,\Omega,\,t}}^2}{||\Delta^{''-1}_{\gamma_{\omega,\,\Omega,\,t}}\bar\partial^\star_{\gamma_{\omega,\,\Omega,\,t}}A_{\omega,\,\Omega,\,t}||_{\gamma_{\omega,\,\Omega,\,t}}}.\end{eqnarray*} Hence, $M(0)= ||\Delta^{''-1}_{\omega}\bar\partial^\star_{\omega}A_{\omega,\,\Omega,\,0}||_{\omega}$. Together with (\ref{eqn:torsion-form_omega-Omega-t_1st-cond_bis_proof_4}), this proves (\ref{eqn:1st-derivative_Gamma_omega_n-1+tOmega}). \hfill $\Box$

\vspace{3ex}

Since $G(\omega_{n-1} + t\Omega) = ||\Gamma_{\omega_{n-1} + t\Omega}||^2_{\gamma_{\omega,\,\Omega,\,t}}$ for every $t\in(-\varepsilon,\,\varepsilon)$ (cf. formula (\ref{eqn:G_energy-functional_balanced})), the above Lemma \ref{Lem:1st-derivative_Gamma_omega_n-1+tOmega} reduces the computation of $\frac{d}{dt}_{|t=0} G(\omega_{n-1} + t\Omega) = \frac{d}{dt}_{|t=0}||\Gamma_{\omega_{n-1} + t\Omega}||^2_{\gamma_{\omega,\,\Omega,\,t}}$ to the computation of $\frac{d}{dt}_{|t=0}||\Gamma_{\omega_{n-1}} + t\,\eta_{\omega,\,\Omega}||^2_{\gamma_{\omega,\,\Omega,\,t}}$. This last derivative at $t=0$ is computed in the following

\begin{Lem}\label{Lem:1st-variation_Gt-tilde} Let $X$ be a compact complex manifold with $\mbox{dim}_\C X=n$. Suppose there exists a balanced metric $\omega$ on $X$.

  For any real $(n-1,\,n-1)$-form $\Omega\in\ker(\bar\partial)$ and for a real $t$ varying in a neighbourhood $(-\varepsilon,\,\varepsilon)$ of $0$ so small that $\omega_{n-1,\,n-1} + t\Omega>0$, if we put $\gamma_t:=\gamma_{\omega,\,\Omega,\,t}$ and $\rho:=\frac{d\gamma_t}{dt}_{|t=0}$, we have: \begin{eqnarray*}\frac{d}{dt}_{|t=0}||\Gamma_{\omega_{n-1}} + t\,\eta_{\omega,\,\Omega}||^2_{\gamma_t} & = & \int\limits_X\eta_{\omega,\,\Omega}\wedge\star_\omega\overline{\Gamma}_{\omega_{n-1}} + \int\limits_X\Gamma_{\omega_{n-1}}\wedge\star_\omega\bigg(\bar\eta_{\omega,\,\Omega} + [\Lambda_\omega,\,\rho\wedge\cdot]\,\overline{\Gamma}_{\omega_{n-1}}\bigg).\end{eqnarray*}

\end{Lem} 

\noindent {\it Proof.} Setting $\alpha_t:=\Gamma_{\omega_{n-1}} + t\,\eta_{\omega,\,\Omega}$, the quantity that we will differentiate at $t=0$ reads: \begin{eqnarray*}||\Gamma_{\omega_{n-1}} + t\,\eta_{\omega,\,\Omega}||^2_{\gamma_t} & = & \int\limits_X|\alpha_t|^2_{\gamma_t}\,dV_{\gamma_t} = \int\limits_X\alpha_t\wedge\star_{\gamma_t}\bar\alpha_t, \hspace{5ex} t\in(-\varepsilon,\,\varepsilon).\end{eqnarray*} Hence \begin{eqnarray}\label{eqn:1st-variation_Gt-tilde_proof_1}\frac{d}{dt}_{|t=0}||\Gamma_{\omega_{n-1}} + t\,\eta_{\omega,\,\Omega}||^2_{\gamma_t} = \int\limits_X\frac{d\alpha_t}{dt}_{|t=0}\wedge\star_{\gamma_0}\bar\alpha_0 + \int\limits_X\alpha_0\wedge\frac{d}{dt}_{|t=0}(\star_{\gamma_t}\bar\alpha_t).\end{eqnarray}

Note that $d\alpha_t/dt = \eta_{\omega,\,\Omega}$. Meanwhile, expanding $\gamma_t$ at $t=0$ and using the fact that $\gamma_0=\omega$ (see (\ref{eqn:root_balanced})), we get $\gamma_t = \omega + t\rho  + O(t^2)$ for $t\in(-\varepsilon,\,\varepsilon)$ if $\varepsilon>0$ is small enough. This implies that $(d/dt)_{|t=0}(\star_{\gamma_t}\bar\alpha_t) = (d/dt)_{|t=0}(\star_{\omega + t\rho}\bar\alpha_t)$. Using Lemma \ref{Lem:1st-variation_Hodge-star} to compute this last derivative at $t=0$, we see that (\ref{eqn:1st-variation_Gt-tilde_proof_1}) translates to: \begin{eqnarray*}\frac{d}{dt}_{|t=0}||\Gamma_{\omega_{n-1}} + t\,\eta_{\omega,\,\Omega}||^2_{\gamma_t} = \int\limits_X\eta_{\omega,\,\Omega}\wedge\star_\omega\overline{\Gamma}_{\omega_{n-1}} + \int\limits_X\Gamma_{\omega_{n-1}}\wedge\star_\omega\bigg(\frac{d\bar\alpha_t}{dt}_{|t=0} + [\Lambda_\omega,\,\rho\wedge\cdot]\,\overline{\Gamma}_{\omega_{n-1}}\bigg),\end{eqnarray*} which proves the contention since $(d\bar\alpha_t/dt)_{|t=0} = \bar\eta_{\omega,\,\Omega}$.  \hfill $\Box$

\vspace{3ex}

Putting Lemmas \ref{Lem:1st-derivative_Gamma_omega_n-1+tOmega} and \ref{Lem:1st-variation_Gt-tilde} together, we get the following

\begin{The}\label{The:1st-variation_Gt-tilde_final} Let $X$ be a compact complex manifold with $\mbox{dim}_\C X=n$. Suppose there exists a balanced metric $\omega$ on $X$.

  For any real $(n-1,\,n-1)$-form $\Omega\in\ker(\bar\partial)$ and for a real $t$ varying in a neighbourhood $(-\varepsilon,\,\varepsilon)$ of $0$ so small that $\omega_{n-1,\,n-1} + t\Omega>0$, we have: \begin{eqnarray}\label{eqn:1st-variation_Gt-tilde_final}\nonumber\frac{d}{dt}_{|t=0} G(\omega_{n-1} + t\Omega)  & = & \langle\langle\eta_{\omega,\,\Omega},\,\Gamma_{\omega_{n-1}}\rangle\rangle_\omega + \langle\langle\Gamma_{\omega_{n-1}},\,\eta_{\omega,\,\Omega} + [\Lambda_\omega,\,\rho\wedge\cdot]\,\Gamma_{\omega_{n-1}}\rangle\rangle_\omega \\
 & + & 2\,||\Gamma_{\omega_{n-1}}||_\omega\,||\Delta_\omega^{''-1}\bar\partial^\star_\omega A_{\omega,\,\Omega,\,0}||_\omega, \end{eqnarray} where $A_{\omega,\,\Omega,\,0}$ is defined by formula (\ref{eqn:A_omega-Omega-0_notation}) and $\rho:=\frac{d\gamma_t}{dt}_{|t=0}$ with $\gamma_t:=\gamma_{\omega,\,\Omega,\,t}$.

In particular, if we choose $\Omega = \omega_{n-1}$, the above formula specialises to \begin{eqnarray}\label{eqn:1st-variation_Gt-tilde_final_special}\frac{d}{dt}_{|t=0} G(\omega_{n-1} + t\omega_{n-1})  = \frac{n+1}{n-1}\,||\Gamma_{\omega_{n-1}}||^2_\omega.\end{eqnarray}

\end{The}

\noindent {\it Proof.} Only formula (\ref{eqn:1st-variation_Gt-tilde_final_special}) still needs a proof. (The definition of the Hodge star operator $\star_\omega$ was used to get the first line of (\ref{eqn:1st-variation_Gt-tilde_final}) from Lemma \ref{Lem:1st-variation_Gt-tilde}.) 

When $\Omega = \omega_{n-1}$, $\eta_{\omega,\,\Omega} = \Gamma_{\omega_{n-1}}$ (see formulae (\ref{eqn:eta_omega-Omega-t_def}) and (\ref{eqn:system_balanced_torsion-form})), while $\rho = \frac{1}{n-1}\,\omega$. This last equality follows from $\rho = \frac{d\gamma_t}{dt}_{|t=0}$ after noting that $(\gamma_t)_{n-1} = (1+t)\,\omega_{n-1}$, hence $\gamma_t = (1+t)^{\frac{1}{n-1}}\,\omega$, when $\Omega = \omega_{n-1}$.

Thus, in this case we get: \begin{eqnarray*}[\Lambda_\omega,\,\rho\wedge\cdot]\,\overline{\Gamma}_{\omega_{n-1}} = \frac{1}{n-1}\,[\Lambda_\omega,\,\omega\wedge\cdot]\,\overline{\Gamma}_{\omega_{n-1}} = -\frac{n-3}{n-1}\,\overline{\Gamma}_{\omega_{n-1}},\end{eqnarray*} since $\overline{\Gamma}_{\omega_{n-1}}$ is of bidegree $(n-2,\,n-1)$ and it is well known that $[\Lambda_\omega,\,\omega\wedge\cdot] = (n-p-q)\,\mbox{Id}$ on $(p,\,q)$-forms (see e.g. [Voi02, Lemma 6.19]). 

Since $||\Gamma_{\omega_{n-1}}||^2_\omega = \int_X\Gamma_{\omega_{n-1}}\wedge\star_\omega\overline{\Gamma}_{\omega_{n-1}}$, (\ref{eqn:1st-variation_Gt-tilde_final_special}) follows at once from these remarks, from (\ref{eqn:1st-variation_Gt-tilde_final}) and from the following observation whose first equality is a consequence of (\ref{eqn:A_omega-Omega-0_notation}): \begin{eqnarray*}A_{\omega,\,\omega_{n-1},\,0} = \bigg(\frac{dP^{(n-1,\,n-1)}_{(1+t)\,\omega_{n-1}}}{dt}_{|t=0}\bigg)(\omega_{n-1}) = 0 \end{eqnarray*} since $P^{(n-1,\,n-1)}_{(1+t)\,\omega_{n-1}}(\omega_{n-1}) =  P^{(n-1,\,n-1)}_{\omega_{n-1}}(\omega_{n-1})$ for every $t$ close enough to $0$. (See (\ref{eqn:P^n-1n-1_omega_n-1_lambda}) for this last fact.)     \hfill $\Box$

\vspace{3ex}

Note that formula (\ref{eqn:1st-variation_Gt-tilde_final_special}) coincides with formula (\ref{eqn:rescaling_G-critical_1+t}) of Lemma \ref{Lem:rescaling_G-critical}. Another conclusion of this study is the following result that was already proved in a different way in (\ref{eqn:critical-points_G-Kaehler}) of Lemma \ref{Lem:rescaling_G-critical}.

\begin{Cor}\label{Cor:crit-points_balanced} Let $X$ be a compact complex $n$-dimensional balanced manifold. 

  Then, for every balanced metric $\omega$ on $X$, $\omega_{n-1}$ is a {\bf critical point} for the functional $G:{\cal B}\longrightarrow[0,\,+\infty)$ if and only if $\omega$ is a {\bf K\"ahler} metric.

\end{Cor}  

\noindent {\it Proof.} If a balanced metric $\omega_{n-1}$ is a critical point for $G$, then $(d_{\omega_{n-1}}G)(\Omega) = \frac{d}{dt}_{|t=0}G(\omega_{n-1} + t\,\Omega) = 0$ for every real $(n-1,\,n-1)$-form $\Omega$ lying in $\ker(\bar\partial)$. In particular, choosing $\Omega:=\omega_{n-1}$ (which is allowed since $\bar\partial\omega_{n-1}=0$), the vanishing of $\frac{d}{dt}_{|t=0}G(\omega_{n-1} + t\,\omega_{n-1})$ implies the vanishing of $\Gamma_{\omega_{n-1}}$ thanks to (\ref{eqn:1st-variation_Gt-tilde_final_special}). On the other hand, the equality $\Gamma_{\omega_{n-1}}=0$ implies that $\omega$ is K\"ahler by Lemma \ref{Lem:Kaehler_balanced_torsion-vanishing}.

Conversely, if $\omega$ is a K\"ahler metric, $\omega_{n-1}$ is a minimiser, hence a critical point, for the functional $G$ by Lemma \ref{Lem:vanishing-G-Kaehler}. 

One can also argue that $\Gamma_{\omega_{n-1}}=0$ if $\omega$ is K\"ahler (again by Lemma \ref{Lem:Kaehler_balanced_torsion-vanishing}), so (\ref{eqn:1st-variation_Gt-tilde_final}) of Theorem \ref{The:1st-variation_Gt-tilde_final} implies that $\frac{d}{dt}_{|t=0}G(\omega_{n-1} + t\,\Omega) = 0$ for every real $(n-1,\,n-1)$-form $\Omega$ lying in $\ker(\bar\partial)$. This means that $\omega_{n-1}$ is a critical point for $G$.  \hfill $\Box$

\vspace{3ex}

In view of a future study of the existence of critical points for the functional $G$, we now give an explicit formula for the $(1,\,1)$-form $\rho=\frac{d\gamma_t}{dt}_{|t=0}$ featuring in (\ref{eqn:1st-variation_Gt-tilde_final}).

\begin{Lem}\label{Lem:balanced_rho_formula}  Let $\omega$ be a Hermitian metric and let $\Omega$ be a real $(n-1,\,n-1)$-form on a compact complex manifold $X$ with $\mbox{dim}_\C X=n$.
  For a real $t$ varying in a neighbourhood $(-\varepsilon,\,\varepsilon)$ of $0$ so small that $\omega_{n-1} + t\Omega>0$, we let $\gamma_t:=\gamma_{\omega,\,\Omega,\,t}$ be the unique positive definite $(1,\,1)$-form on $X$ such that $(\gamma_t)_{n-1} = \omega_{n-1} + t\Omega$ and we set $\rho:=\frac{d\gamma_t}{dt}_{|t=0}$.

Then, the following formula holds: \begin{eqnarray}\label{eqn:balanced_rho_formula}\rho = \frac{1}{n-1}\,\Lambda_\omega(\star_\omega\Omega)\,\omega - \star_\omega\Omega.\end{eqnarray}

\end{Lem} 

\noindent {\it Proof.} Since the result is of a pointwise nature, we fix an arbitrary point $x\in X$ and local coordinates $z_1,\dots , z_n$ centred at $x$ such that \begin{eqnarray*}\omega(x) = \sum\limits_{j=1}^n idz_j\wedge d\bar{z}_j  \hspace{3ex} \mbox{and} \hspace{3ex} \Omega(x) = \sum\limits_{j=1}^n\Omega_j\,\widehat{idz_j\wedge d\bar{z}_j},\end{eqnarray*} where, for every $j$, $\Omega_j\in\R$ and $\widehat{idz_j\wedge d\bar{z}_j}$ is the $(n-1,\,n-1)$-form given as the product of all the $(1,\,1)$-forms $idz_k\wedge d\bar{z}_k$ with $k\neq j$. Then \begin{eqnarray*}(\omega_{n-1} + t\Omega)(x) = \sum\limits_{j=1}^n(1 + t\Omega_j)\,\widehat{idz_j\wedge d\bar{z}_j}, \hspace{5ex} t\in(-\varepsilon,\,\varepsilon),\end{eqnarray*} and the $(n-1)$-st root of this form is given by \begin{eqnarray*}\gamma_t(x) = \sum\limits_{j=1}^n\frac{\Pi_{k=1}^n(1 + t\Omega_k)^{\frac{1}{n-1}}}{1 + t\Omega_j}\,idz_j\wedge d\bar{z}_j, \hspace{5ex} t\in(-\varepsilon,\,\varepsilon).\end{eqnarray*}

Differentiating with respect to $t$ and then evaluating at $t=0$, we get: 

\begin{eqnarray*}\frac{d\gamma_t}{dt}_{|t=0}(x) & = & - \sum\limits_{j=1}^n\frac{\Omega_j\,\Pi_{k=1}^n(1 + t\Omega_k)^{\frac{1}{n-1}}}{(1 + t\Omega_j)^2}_{\bigg|t=0}\,idz_j\wedge d\bar{z}_j \\
 & + & \frac{1}{n-1}\,\sum\limits_{j=1}^n\frac{\sum\limits_{k=1}^n(1 + t\Omega_1)^{\frac{1}{n-1}}\cdots\widehat{(1 + t\Omega_k)^{\frac{1}{n-1}}}\cdots(1 + t\Omega_n)^{\frac{1}{n-1}}\,(1 + t\Omega_k)^{\frac{2-n}{n-1}}\,\Omega_k}{1 + t\Omega_j}_{\bigg|t=0}\,idz_j\wedge d\bar{z}_j \\
 & = & -\sum\limits_{j=1}^n\Omega_j\,idz_j\wedge d\bar{z}_j + \frac{1}{n-1}\,\sum\limits_{j=1}^n\bigg(\sum\limits_{k=1}^n\Omega_k\bigg)\,idz_j\wedge d\bar{z}_j \\
 & = & \frac{1}{n-1}\,\bigg(\sum\limits_{k=1}^n\Omega_k\bigg)\,\omega(x) - \sum\limits_{j=1}^n\Omega_j\,idz_j\wedge d\bar{z}_j.\end{eqnarray*} This proves (\ref{eqn:balanced_rho_formula}) since $(\star_\omega\Omega)(x) = \sum_{j=1}^n\Omega_j\,idz_j\wedge d\bar{z}_j$ and $\Lambda_\omega(\star_\omega\Omega)(x) = \sum_{k=1}^n\Omega_k$.  \hfill $\Box$

\section{Simultaneous study of SKT and balanced metrics}\label{section:simultaneous_SKT-bal}

Having separately studied SKT and balanced metrics in the previous sections, we indicate in this short section a possible simultaneous variational approach to these two classes of metrics in vew of a future attack on the problem proposed by Fino and Vezzoni in [FV15] asking whether or not the existence of both kinds of metrics on a compact complex manifold implies the existence of a K\"ahler metric.

\subsection{An SKT-balanced observation}\label{subsection:SKT-bal_obs} We now run again the argument of [Pop15, Proposition 1.1], made even shorter in [DP20, Proposition 2.6.], but with two different metrics. The next observation will be the starting point of the subsequent considerations.

\begin{Prop}\label{Prop:SKT-bal_obs} Let $X$ be a compact complex manifold with $\mbox{dim}_\C X =n$. For any {\bf SKT} metric $\omega$ (if any) and any {\bf balanced} metric $\gamma$ (if any) on $X$, we have: \begin{equation}\label{eqn:SKT-bal_obs}\langle\langle\bar\partial\omega,\,\bar\partial\gamma\rangle\rangle_\gamma = 0,\end{equation} where $\langle\langle\cdot\,\,,\,\,\cdot\rangle\rangle_\gamma$ is the $L^2$ inner product induced by the metric $\gamma$. 

\end{Prop}

\noindent {\it Proof.} The {\it SKT} assumption on $\omega$ translates to any of the following equivalent properties: \begin{eqnarray}\label{eqn:pluriclosed-equiv}\partial\bar\partial\omega=0 \Longleftrightarrow \partial\omega\in\ker\bar\partial \Longleftrightarrow \star_\omega(\partial\omega)\in\ker\partial^{\star}_\omega,\end{eqnarray}

\noindent where the last equivalence follows from the standard formula $\partial^{\star}_\omega=-\star_\omega\bar\partial\star_\omega$ involving the Hodge-star isomorphism $\star_{\omega}:\Lambda^{p,\,q}T^{\star}X\rightarrow \Lambda^{n-q,\,n-p}T^{\star}X$ defined by $\omega$ for arbitrary $p,q=0,\dots , n$.

 Meanwhile, the {\it balanced} assumption on $\gamma$ translates to any of the following equivalent properties: \begin{eqnarray*}\label{eqn:bal-equiv}d\gamma^{n-1}=0 \Longleftrightarrow \partial\gamma^{n-1}=0 \Longleftrightarrow \gamma^{n-2}\wedge\partial\gamma = 0   \Longleftrightarrow \partial\gamma \,\,\mbox{is $\gamma$-primitive}.\end{eqnarray*}

Now, recall the following standard formula (cf. e.g. [Voi02, Proposition 6.29, p. 150]) for the Hodge star operator $\star = \star_\omega$ of any Hermitian metric $\omega$ applied to {\it primitive} forms $v$ of arbitrary bidegree $(p, \, q)$: \begin{eqnarray}\label{eqn:prim-form-star-formula-gen}\star\, v = (-1)^{k(k+1)/2}\, i^{p-q}\, \frac{\omega^{n-p-q}\wedge v}{(n-p-q)!}, \hspace{2ex} \mbox{where}\,\, k:=p+q.\end{eqnarray}

 Thus, since the $(2,\,1)$-form $\partial\gamma$ is $\gamma$-primitive when $\gamma$ is balanced, the standard formula (\ref{eqn:prim-form-star-formula-gen}) yields: \begin{eqnarray}\label{eqn:consequence_balanced}\star_\gamma(\partial\gamma) = i\,\frac{\gamma^{n-3}}{(n-3)!}\wedge\partial\gamma = \frac{i}{(n-2)!}\,\partial\gamma^{n-2}\in\mbox{Im}\,\partial.\end{eqnarray}

\vspace{1ex}

From the subspaces $\ker\partial^{\star}_\omega$ and $\mbox{Im}\,\partial$ of $C^\infty_{n-1,\,n-2}(X,\,\C)$ being $L^2_\omega$-orthogonal and from (\ref{eqn:pluriclosed-equiv}) and (\ref{eqn:consequence_balanced}) we infer the first relation below: \begin{eqnarray*}\star_\omega(\partial\omega)\perp_{L^2_\omega}\star_\gamma(\partial\gamma) & \iff & \langle\langle\star_\gamma(\partial\gamma),\,\star_\omega(\partial\omega)\rangle\rangle_\omega = 0 \iff \int\limits_X\bigg\langle\star_\gamma(\partial\gamma),\,\star_\omega(\partial\omega)\bigg\rangle_\omega\,dV_\omega = 0 \\
  & \iff & \int\limits_X\star_\gamma(\partial\gamma)\wedge\star_\omega\star_\omega(\bar\partial\omega) = 0  \iff \int\limits_X\bar\partial\omega\wedge\star_\gamma\overline{(\bar\partial\gamma)} = 0 \\
  & \iff & \int\limits_X\bigg\langle\bar\partial\omega,\,\bar\partial\gamma\bigg\rangle_\gamma\,dV_\gamma = 0 \iff \langle\langle\bar\partial\omega,\,\bar\partial\gamma\rangle\rangle_\gamma = 0.\end{eqnarray*}

This proves the contention.  \hfill $\Box$

\hspace{2ex}

Note that, if $\omega=\gamma$, we deduce from Proposition \ref{Prop:SKT-bal_obs} the well-known fact that, if a Hermitian metric $\omega$ is at once SKT and balanced, it is K\"ahler. We even get the following slight generalisation of this observation.

\begin{Cor}\label{Cor:SKT-bal_obs} Let $X$ be a compact complex manifold with $\mbox{dim}_\C X =n$. Suppose there exists an {\bf SKT} metric $\omega$ and a {\bf balanced} metric $\gamma$ on $X$ such that $$\bar\partial\omega = \bar\partial\gamma.$$ Then, both $\omega$ and $\gamma$ are {\bf K\"ahler}.

\end{Cor}  

\noindent {\it Proof.} If $\bar\partial\omega = \bar\partial\gamma$, (\ref{eqn:SKT-bal_obs}) reads $||\bar\partial\omega||^2_\gamma = ||\bar\partial\gamma||^2_\gamma=0$. Hence, $\bar\partial\omega = \bar\partial\gamma =0$, translating the fact that $\omega$ and $\gamma$ are K\"ahler.  \hfill $\Box$

\vspace{2ex}

It also follows from Proposition \ref{Prop:SKT-bal_obs} that, whenever $\omega$ is an {\it SKT metric} and $\gamma$ is a {\it balanced metric}, we have: \begin{equation}\label{eqn:SKT-bal_difference-orth}||\bar\partial\omega - \bar\partial\gamma||^2_\gamma = ||\bar\partial\omega||^2_\gamma + ||\bar\partial\gamma||^2_\gamma .\end{equation}



\subsection{The SKT-balanced functional}\label{subsection:functional_skt-balanced}

We now introduce a functional depending on both SKT and balanced metrics that are assumed to exist on a given compact complex manifold.



\begin{Def}\label{Def:functional_skt-balanced} Let $X$ be a compact complex manifold with $\mbox{dim}_\C X=n$. Suppose both {\bf SKT} and {\bf balanced} metrics
 exist on $X$. For a SKT metric $\omega$ and a $(\gamma)_{n-1}$ of any balanced metric $\gamma$, we define the following {\bf energy functional}:
  
 \begin{equation}\label{eqn:H_energy-functional_skt-balanced}H : {\cal S}\times{\cal B} \to [0,\,+\infty), \hspace{3ex} H(\omega,\,\gamma_{n-1}) = ||\Lambda_\omega(\partial\omega)||^2_\gamma,\end{equation}

\noindent where $||\,\,\,||_\gamma$ is the $L^2$-norm induced by $\gamma$. 

\end{Def}

We have: \begin{equation}\label{eqn:H_energy-functional_skt-balanced_bis}H(\omega,\,\gamma_{n-1}) = \int\limits_X|\Lambda_\omega(\partial\omega)|^2_\gamma\,dV_\gamma = \int\limits_X\Lambda_\omega(\partial\omega)\wedge\star_\gamma\Lambda_\omega(\bar\partial\omega) = i\int\limits_X\Lambda_\omega(\partial\omega)\wedge\Lambda_\omega(\bar\partial\omega)\wedge\gamma_{n-1}\end{equation} for every SKT metric $\omega$ and every balanced metric $\gamma$.

\vspace{2ex}

The first trivial observation that justifies the introduction of the functional $H$ is the following.

\begin{Lem}\label{Lem:vanishing-H-omega-Kaehler} Let $\omega_0$ be an SKT metric and let $\gamma_0$ be a balanced metric on $X$. For any $\omega\in{\cal S}$ the following equivalence holds: \begin{equation}\label{eqn:vanishing-G-omega-Kaehler}\omega \hspace{1ex} \mbox{is K\"ahler} \iff H(\omega,\,(\gamma_0)_{n-1})=0.\end{equation}

\end{Lem} 

\noindent {\it Proof.} The condition $H(\omega,\,(\gamma_0)_{n-1})=0$ is equivalent to $\Lambda_\omega(\partial\omega)=0$, which in turn is equivalent to $\partial\omega$ being primitive w.r.t. $\omega$. This last property translates to $\omega^{n-2}\wedge\partial\omega = 0$ (since $\partial\omega$ is a $3$-form), which amounts to $\partial\omega^{n-1} = 0$. This means that $\omega$ satisfies the balanced condition. Being already SKT, the metric $\omega$ is balanced if and only if it is K\"ahler. (See Corollary \ref{Cor:SKT-bal_obs}.) \hfill $\Box$

\vspace{3ex}

The idea is to use the balanced metric $\gamma_{n-1}$ as a reference metric and to vary the SKT metric $\omega$. Note that, for any Hermitian metric $\omega$ and every real $\lambda>0$, we have $\Lambda_{\lambda\omega} = (1/\lambda)\,\Lambda_\omega$ as operators acting in any degree, hence \begin{eqnarray*}H(\lambda\omega,\,(\gamma_0)_{n-1}) = H(\omega,\,(\gamma_0)_{n-1}), \hspace{5ex} \lambda>0,\end{eqnarray*} showing that the functional $H$ is constant in its first variable along the rays of the cone ${\cal S}$ whenever its second variable $\gamma_{n-1}$ has been fixed.

We now to compute the first variation of $H$ in the $\omega$ variable. We fix an SKT metric $\omega$ on $X$ and we vary it in $\cal S$ along the path $\omega + t\eta$, where $\eta\in C^\infty_{1,\,1}(X,\,\R)$ is a fixed real $(1,\,1)$-form. We will be considering reals $t\in(-\varepsilon,\,\varepsilon)$ with $\varepsilon>0$ small enough to ensure that $\omega + t\eta>0$ (i.e. $\omega + t\eta$ is a, necessarily SKT, metric on $X$).

Using (\ref{eqn:H_energy-functional_skt-balanced_bis}), we get: \begin{eqnarray}\label{eqn:H_energy-functional_skt-balanced_computation_1}\nonumber\frac{d}{dt}\bigg|_{t=0} H(\omega + t\eta,\,\gamma_{n-1}) & = & i\,\frac{d}{dt}\bigg|_{t=0}\int\limits_X\Lambda_{\omega + t\eta}(\partial\omega + t\,\partial\eta)\wedge\Lambda_{\omega + t\eta}(\bar\partial\omega + t\,\bar\partial\eta)\wedge\gamma_{n-1}  \\
  \nonumber & = & i\int\limits_X\frac{d}{dt}\bigg|_{t=0}\bigg(\Lambda_{\omega + t\eta}(\partial\omega + t\,\partial\eta)\bigg)\wedge\Lambda_\omega(\bar\partial\omega)\wedge\gamma_{n-1} \\
  \nonumber & + & i\int\limits_X\Lambda_\omega(\partial\omega)\wedge\frac{d}{dt}\bigg|_{t=0}\bigg(\Lambda_{\omega + t\eta}(\bar\partial\omega + t\,\bar\partial\eta)\bigg)\wedge\gamma_{n-1}.\end{eqnarray}

We now apply Lemma \ref{Lem:1st-variation_trace_pq} to the computation of the derivatives at $t=0$ of the two factors of the shape $\Lambda_{\omega + t\eta}(\cdot)$. Formula (\ref{eqn:1st-variation_trace_pq}) yields: \begin{eqnarray*}\frac{d}{dt}\bigg|_{t=0}\bigg(\Lambda_{\omega + t\eta}(\partial\omega + t\,\partial\eta)\bigg) & = & \Lambda_\omega(\partial\eta) - (\eta\wedge\cdot)^\star_\omega(\partial\omega) \\
 \frac{d}{dt}\bigg|_{t=0}\bigg(\Lambda_{\omega + t\eta}(\bar\partial\omega + t\,\bar\partial{\eta})\bigg) & = & \Lambda_\omega(\bar\partial\eta) - (\eta\wedge\cdot)^\star_\omega(\bar\partial\omega).\end{eqnarray*}

\begin{Lem}\label{Lem:1st-variation_H_omega} Let $X$ be a compact complex manifold with $\mbox{dim}_\C X=n$. Suppose there exists an SKT metric $\omega$ and a balanced metric $\gamma$ on $X$.

  For any real $(1,\,1)$-form  $\eta\in C^\infty_{1,\,1}(X,\,\R)$, and for a real $t$ varying in a small neighbourhood of $0$ such that $\omega + t\eta>0$, we have \begin{eqnarray}\label{eqn:1st-variation_G_omega}\nonumber\frac{d}{dt}\bigg|_{t=0} H(\omega + t\eta,\,\gamma_{n-1}) & = & 2\,\mbox{Re}\,\bigg(i\int\limits_X\Lambda_\omega(\partial\eta)\wedge\Lambda_\omega(\bar\partial\omega)\wedge\gamma_{n-1}\bigg) \\
    & - & 2\,\mbox{Re}\,\bigg(i\int\limits_X(\eta\wedge\cdot)^\star_\omega(\partial\omega)\wedge\Lambda_\omega(\bar\partial\omega)\wedge\gamma_{n-1}\bigg).\end{eqnarray}

\end{Lem}

\vspace{3ex}

Note that, if we take $\eta=\omega$, we get $(\eta\wedge\cdot)^\star_\omega = \Lambda_\omega$, so the right-hand side term of (\ref{eqn:1st-variation_G_omega}) vanishes. This shows that $(d/dt)_{t=0}H((1+t)\omega,\,\gamma_{n-1}) = 0$, in line with the above observation that $H$ is constant in the first variable along the rays $\{\lambda\omega\,\mid\,\lambda>0\}$ of the cone ${\cal S}$ for any fixed second variable $\gamma_{n-1}$.

\subsection{Normalised SKT functionals}\label{subsection:normalised_SKT-functionals} Another natural variational approach to the Fino-Vezzoni conjecture [FV15] (predicting that the simultaneous existence of an SKT metric and of a balanced metric on a compact complex manifold ought to imply the existence of a K\"ahler metric) seems to be via an appropriate normalisation of our SKT functional defined in $\S$\ref{subsection:functional_skt-Delta}.

\begin{Def}\label{Def:normalised_SKT-functionals} Let $X$ be a compact complex SKT manifold with $\mbox{dim}_\C X=n$. Fix a Hermitian metric $\nu$ on $X$. We define the following $\nu$-dependent functional acting on the SKT metrics of $X$: \begin{eqnarray}\label{eqn:normalised_SKT-functional}\widetilde{F}_{\nu}: {\cal S}\rightarrow[0,\,+\infty), \hspace{5ex} \widetilde{F}_{\nu}(\omega):= \frac{F(\omega)}{\bigg(\int_{X}\omega\wedge\nu_{n-1}\bigg)^n},\end{eqnarray} where $F$ is the functional introduced in Definition \ref{Def:F_energy-functional_SKT}.

\end{Def}

The behaviour of $F$ under rescalings observed in (\ref{eqn:rescaling_F-critical}) implies the scaling invariance of $\widetilde{F}_{\nu}$: \begin{eqnarray*}\widetilde{F}_{\nu}(\lambda\omega) = \widetilde{F}_{\nu}(\omega)\end{eqnarray*} for every real $\lambda>0$. On the other hand, the background metric $\nu$ can be chosen to be balanced (assuming such metrics exist on $X$) and even made to range over a certain family of metrics. We will hopefully investigate these aspects in future work.

\vspace{2ex}

We now observe a formula for the first variation of the normalised functional $\widetilde{F}_{\nu}$ in terms of the analogous expression for the unnormalised functional $F$ that was computed in Theorem \ref{The:1st-variation_Ft-tilde_final}.

\begin{Prop}\label{Prop:1st-variation-normalised_SKT-functional} Let $X$ be a compact complex manifold with $\mbox{dim}_\C X=n$. Fix a Hermitian metric $\nu$ on $X$. Then, for any SKT metric $\omega$ (if any) and any real $(1,\,1)$-form $\gamma\in\ker(\partial\bar\partial)$ on $X$, we have:
  \begin{eqnarray}\label{eqn:1st-variation-normalised_SKT-functional}(d_{\omega}\widetilde{F}_{\nu})(\gamma) = \frac{1}{\bigg(\int_{X}\omega\wedge\nu_{n-1}\bigg)^n}\,\bigg((d_{\omega}F)(\gamma) - n\,\frac{\int_{X}\gamma\wedge\nu_{n-1}}{\int_{X}\omega\wedge\nu_{n-1}}\,F(\omega)\bigg),\end{eqnarray} where $(d_{\omega}F)(\gamma)$ is given by formula (\ref{eqn:1st-variation_Ft-tilde_final}) in Theorem \ref{The:1st-variation_Ft-tilde_final}.

\end{Prop}

\noindent {\it Proof.} Straightforward computations yield: \begin{eqnarray*}(d_{\omega}\widetilde{F}_{\nu})(\gamma) & = & \frac{d}{dt}\bigg[\frac{1}{\bigg(\int_{X}(\omega + t\gamma)\wedge\nu_{n-1}\bigg)^n}\,F(\omega + t\gamma)\bigg]_{t=0}  = \frac{1}{\bigg(\int_{X}\omega\wedge\nu_{n-1}\bigg)^n}\,(d_{\omega}F)(\gamma) \\
  & & - \frac{1}{\bigg(\int_{X}\omega\wedge\nu_{n-1}\bigg)^{2n}}\,n\,\bigg(\int_{X}\omega\wedge\nu_{n-1}\bigg)^{n-1}\,\bigg(\int_{X}\gamma\wedge\nu_{n-1}\bigg)\,F(\omega).\end{eqnarray*} This is formula (\ref{eqn:1st-variation-normalised_SKT-functional}).  \hfill $\Box$

\vspace{3ex}

The normalised functionals $\widetilde{F}_{\nu}$ may lead to singling out obstructions to the K\"ahlerianity of a compact SKT manifold $X$ in the following way. Let $n$ be the complex dimension of $X$ and let $\nu$ be a Hermitian metric thereon. We restrict attention to the set $V_\nu$ of SKT metrics $\omega$ on $X$ that are normalised in the following way with respect to $\nu$: $$\int\limits_X\omega\wedge\nu_{n-1} = 1.$$  

Then, $\widetilde{F}_\nu(\omega) = F(\omega)$ for every $\omega\in V_\nu$ (see Definition \ref{Def:normalised_SKT-functionals}). Now, the scaling-invaraint functional $\widetilde{F}_\nu$ is determined by its restriction to $V_\nu$. We put: $$c_\nu:=\inf\limits_{\omega\in{\cal S}}\widetilde{F}_\nu(\omega) = \inf\limits_{\omega\in V_\nu}\widetilde{F}_\nu(\omega) = \inf\limits_{\omega\in V_\nu}F(\omega)\geq 0.$$

For every $\varepsilon>0$, let $\omega_\varepsilon\in V_\nu$ such that $c_\nu\leq \widetilde{F}_\nu(\omega_\varepsilon)< c_\nu + \varepsilon$. Since the set $V_\nu$ is relatively compact in the space of positive (see terminology in e.g. [Dem97, III-1.B.]) $(1,\,1)$-currents with respect to the weak topology of currents, there exists a subsequence $\varepsilon_k\downarrow 0$ and a positive $(1,\,1)$-current $T_\nu\geq 0$ on $X$ such that the sequence $(\omega_{\varepsilon_k})_k$ converges weakly to $T_\nu$ as $k\to +\infty$. By construction, we have: $$\int\limits_XT_\nu\wedge\nu_{n-1} = 1$$ and $\partial\bar\partial T_\nu=0$ (since $\partial\bar\partial\omega_\varepsilon=0$ for every $\varepsilon$). 

The possible failure of the current $T_\nu\geq 0$ to be either a $C^\infty$ form or strictly positive (for example, in the sense that it is bounded below by a positive multiple of a Hermitian metric on $X$, in which case we can call it ``an SKT current minimising $\widetilde{F}_\nu$'') constitutes an obstruction to the existence of minimisers for the functional $\widetilde{F}_\nu$. Should it turn out that the critical points of $\widetilde{F}_\nu$, if any, are precisely the K\"ahler metrics of $X$, if any, these critical points will further coincide with the minimisers of $\widetilde{F}_\nu$. In that case, the currents $T_\nu$ can be seen as obstructions to the existence of K\"ahler metrics on $X$.

\vspace{3ex}

\noindent {\bf References.} \\

\vspace{1ex}

\noindent [AB95]\, L. Alessandrini, G. Bassanelli --- {\it Modifications of compact balanced manifolds} --- C. R. Acad. Sci. Paris {\bf 320} (1995), 1517-1522.

\vspace{1ex}

\noindent [AI01]\, B. Alexandrov, S. Ivanov --- {\it Vanishing Theorems on Hermitian Manifolds} --- Differential Geom. Appl. {\bf 14} (2001), no. 3, 251–265.

\vspace{1ex}

\noindent [AN22]\, R.M. Arroyo, M. Nicolini --- {\it SKT Structures on Nilmanifolds} --- Math. Z. {\bf 302} (2022), 1307-1320.

\vspace{1ex}

\noindent [BP18]\, H. Bellitir, D. Popovici --- {\it Positivity Cones under Deformations of Complex Structures} --- Riv. Mat. Univ. Parma, Vol. {\bf 9} (2018), 133-176.

\vspace{1ex}

\noindent [Chi14]\, I. Chiose --- {\it Obstructions to the Existence of K\"ahler Structures on Compact Complex Manifolds} --- Proc. Amer. Math. Soc. {\bf 142} (2014), no. 10, 3561–3568.

\vspace{1ex}

\noindent [Dem97]\, J.-P. Demailly --- {\it Complex Analytic and Algebraic Geometry} --- \url{https://www-fourier.ujf-grenoble.fr/~demailly/manuscripts/agbook.pdf}

\vspace{1ex}

\noindent [DP20]\, S. Dinew, D. Popovici --- {\it A Generalised Volume Invariant for Aeppli Cohomology Classes of Hermitian-Symplectic Metrics} ---  Adv. Math. {\bf 393} (2021), Paper No. 108056, 46 pp.

\vspace{1ex}

\noindent[FPS04] A. Fino, M. Parton, S. Salamon --- {\it Families of Strong KT Structures in Six Dimensions}--- Comment. Math. Helv. {\bf 79} (2004), 317-340.

\vspace{1ex}

\noindent [FV15]\, A. Fino, L. Vezzoni --- {\it Special Hermitian metrics on compact solvmanifolds} --- J. Geom. Phys. {\bf 91} (2015), 40-53.

\vspace{1ex}

\noindent [FV16]\, A. Fino, L. Vezzoni --- {\it On the Existence of Balanced and SKT Metrics on Nilmanifolds} --- Proc. Amer. Math. Soc. {\bf 144} (2016), no. 6, 2455-2459. 

\vspace{1ex}

\noindent [FLY12]\, J.Fu, J.Li, S.-T. Yau --- {\it Balanced Metrics on Non-K\"ahler Calabi-Yau Threefolds} --- J. Diff. Geom. {\bf 90} (2012) 81-129.

\vspace{1ex}

\noindent [Gau77a]\, P. Gauduchon --- {\it Le th\'eor\`eme de l'excentricit\'e nulle} --- C. R. Acad. Sci. Paris, S\'er. A, {\bf 285} (1977), 387-390.

\vspace{1ex}

\noindent [Gau77b]\, P. Gauduchon --- {\it Fibr\'es hermitiens \`a endomorphisme de Ricci non n\'egatif} --- Bull. Soc. Math. France {\bf 105} (1977) 113-140.

\vspace{1ex}

\noindent [IP12]\, S. Ivanov, G. Papadopoulos --- {\it Vanishing Theorems on $(l/k)$-strong K\"ahler Manifolds with Torsion} --- arXiv e-print DG 1202.6470v1.

\vspace{1ex}

\noindent [Kod86]\, K. Kodaira --- {\it Complex Manifolds and Deformations of Complex Structures} --- Grundlehren der Math. Wiss. {\bf 283}, Springer (1986).

\vspace{1ex}

\noindent [KS60]\, K. Kodaira, D.C. Spencer --- {\it On Deformations of Complex Analytic Structures, III. Stability Theorems for Complex Structures} --- Ann. of Math. {\bf 71}, no.1 (1960), 43-76.

\vspace{1ex}

\noindent [MP21]\, {\it Balanced Hyperbolic and Divisorially Hyperbolic Compact Complex Manifolds} --- arXiv e-print CV 2107.08972v2, to appear in Mathematical Research Letters.

\vspace{1ex}

\noindent [Mic83]\, M. L. Michelsohn --- {\it On the Existence of Special Metrics in Complex Geometry} --- Acta Math. {\bf 143} (1983) 261-295.

\vspace{1ex}

\noindent [Pop13]\, D. Popovici --- {\it Deformation Limits of Projective Manifolds: Hodge Numbers and Strongly Gauduchon Metrics} --- Invent. Math. {\bf 194} (2013), 515-534.

\vspace{1ex}

\noindent [Pop15]\, D. Popovici --- {\it Aeppli Cohomology Classes Associated with Gauduchon Metrics on Compact Complex Manifolds} --- Bull. Soc. Math. France {\bf 143}, no. 3 (2015), 1-37.

\vspace{1ex}

\noindent [Pop22]\, D. Popovici --- {\it Pluriclosed Star Split Hermitian Metrics} --- arXiv:2211.10267v1 

\vspace{1ex}

\noindent [Sch07]\, M. Schweitzer --- {\it Autour de la cohomologie de Bott-Chern} --- arXiv e-print math.AG/0709.3528v1.

\vspace{1ex}

\noindent[ST10]\, J. Streets, G. Tian --- {\it Parabolic flow of pluriclosed metrics}--- Int. Math. Res. Not. IMRN {\bf 16} (2010), 3101-3133.

\vspace{1ex}

\noindent [Ver14]\, M. Verbitsky --- {\it Rational Curves and Special Metrics on Twistor Spaces} --- Geometry and Topology {\bf 18} (2014), 897–909.

\vspace{1ex}

\noindent [Voi02]\, C. Voisin --- {\it Hodge Theory and Complex Algebraic Geometry. I.} --- Cambridge Studies in Advanced Mathematics, 76, Cambridge University Press, Cambridge, 2002.

\vspace{1ex}

\vspace{6ex}

\noindent Department of Mathematics and Computer Science     \hfill Institut de Math\'ematiques de Toulouse,

\noindent Jagiellonian University     \hfill  Universit\'e Paul Sabatier,

\noindent 30-409 Krak\'ow, Ul. Lojasiewicza 6, Poland    \hfill  118 route de Narbonne, 31062 Toulouse, France

\noindent  Email: Slawomir.Dinew@im.uj.edu.pl                \hfill     Email: popovici@math.univ-toulouse.fr

\end{document}